\pdfoutput=1
\RequirePackage{fix-cm}
\documentclass{amsart}
\usepackage{geometry}
\usepackage{microtype}\microtypesetup{final,verbose=errors}
\usepackage{inputenc}\inputencoding{utf8}
\usepackage{comment}
\newcommand{\NEWOLD}[2]{%
	\expandafter\newcommand\csname NEW#1\endcsname[1]{\csname #1\endcsname=##1}%
	\expandafter\newcommand\csname OLD#1\endcsname{\csname #1\endcsname=#2}%
	}
\NEWOLD{pretolerance}{100}
\NEWOLD{tolerance}{200}
\NEWOLD{hbadness}{1000}
\NEWOLD{vbadness}{1000}
\AtBeginDocument{%
	\NEWpretolerance{-1}%
	\NEWtolerance{100}%
	\NEWhbadness{100}%
	\NEWvbadness{100}%
	}
\NEWOLD{fboxsep}{3pt}
\NEWOLD{fboxrule}{.4pt}
\NEWOLD{arraycolsep}{5pt}
\NEWOLD{tabcolsep}{6pt}
\AtBeginDocument{%
	\NEWfboxrule{0.004pt}%
	\NEWfboxsep{-\fboxrule}%
	\NEWarraycolsep{0pt}%
	\NEWtabcolsep{0pt}%
	}
\NEWOLD{leftmargini}{2.5em}
\NEWOLD{leftmarginii}{2.2em}
\NEWOLD{leftmarginiii}{1.87em}
\NEWOLD{leftmarginiv}{1.7em}
\AtBeginDocument{%
	\NEWleftmargini{12.7778pt}%
	\NEWleftmarginii{12.7778pt}%
	\NEWleftmarginiii{18.61116pt}%
	\NEWleftmarginiv{15.62502pt}%
	}
\NEWOLD{abovedisplayskip}{10pt plus 2pt minus 5pt}
\NEWOLD{belowdisplayskip}{\abovedisplayskip}
\NEWOLD{abovedisplayshortskip}{0pt plus 3pt}
\NEWOLD{belowdisplayshortskip}{6.5pt plus 3.5pt minus 3pt}
\AtBeginDocument{%
	\NEWabovedisplayshortskip{0pt plus 2pt minus 5pt}%
	\NEWbelowdisplayshortskip{10pt plus 2pt minus 5pt}%
	}
\makeatletter
\let\fboxT\fbox
\newcommand{\fboxM}[1]{\mathpalette\fboxM@{#1}}
\newcommand{\fboxM@}[2]{\fboxT{$\mathsurround=0pt#1#2$}}
\renewcommand{\fbox}[1]{\ifmmode\fboxM{#1}\else\fboxT{#1}\fi}
\newcommand{\swapabovedisplayskip}{\SwapAboveDisplaySkip}
\makeatother
\usepackage{mathtools}
\DeclarePairedDelimiter{\sM}{\{}{\}}
\DeclarePairedDelimiter{\pM}{(}{)}
\DeclarePairedDelimiter{\bM}{[}{]}
\DeclarePairedDelimiter{\lbrpM}{[}{)}
\DeclarePairedDelimiter{\abM}{\langle}{\rangle}
\DeclarePairedDelimiter{\vM}{\lvert}{\rvert}
\DeclarePairedDelimiter{\vvM}{\lVert}{\rVert}
\let\lim\undefined
\let\sup\undefined
\let\inf\undefined

\let\max\undefined
\let\min\undefined

\let\cos\undefined

\let\log\undefined

\DeclareMathOperator*{\lim}{lim}
\DeclareMathOperator*{\sup}{sup}
\DeclareMathOperator*{\inf}{inf}

\DeclareMathOperator*{\max}{max}
\DeclareMathOperator*{\min}{min}

\DeclareMathOperator{\cos}{cos}

\DeclareMathOperator{\log}{log}

\DeclareMathOperator{\supp}{supp}
\DeclareMathOperator{\spec}{spec}
\DeclareMathOperator{\dist}{dist}

\usepackage{amssymb}
\DeclareMathAlphabet{\mathbbm}{U}{bbm}{m}{n}
\SetMathAlphabet\mathbbm{bold}{U}{bbm}{bx}{n}
\DeclareMathAlphabet{\mathbbmss}{U}{bbmss}{m}{n}
\SetMathAlphabet\mathbbmss{bold}{U}{bbmss}{bx}{n}
\DeclareSymbolFont{rsfs}{U}{rsfs}{m}{n}
\DeclareSymbolFontAlphabet{\mathscr}{rsfs}
\usepackage{mleftright}

\let\left\mleft
\let\right\mright
\usepackage{tensor}

\usepackage{bm}
\usepackage{xcolor}
\usepackage{accsupp}
\usepackage{pagecolor}
\colorlet{Bc}{\thepagecolor}
\colorlet{Tc}{.}
\let\originalcolorbox\colorbox
\let\originalfcolorbox\fcolorbox
\renewcommand{\colorbox}[2]{\colorlet{Bc}{#1}\originalcolorbox{#1}{#2}}
\renewcommand{\fcolorbox}[3]{\colorlet{Bc}{#2}\originalfcolorbox{#1}{#2}{#3}}
\makeatletter
\renewcommand{\,}{\@ifstar\commaStar\commaNoStar}
\newcommand{\commaStar}{\mathchoice
	{,\sbox0{$\mathsurround=0pt\textstyle,$}\kern-3.6\ht0\rlap{\color{Bc}\smash{\rule[-0.999\dp0]{1.4\ht0}{1.1\ht0+1.1\dp0}}}\kern2.021\ht0\llap{\BeginAccSupp{ActualText=}\textup{,}\EndAccSupp{}}}
	{,\sbox0{$\mathsurround=0pt\textstyle,$}\kern-3.6\ht0\rlap{\color{Bc}\smash{\rule[-0.999\dp0]{1.4\ht0}{1.1\ht0+1.1\dp0}}}\kern2.021\ht0\llap{\BeginAccSupp{ActualText=}\textup{,}\EndAccSupp{}}}
	{\textup{,}}
	{\textup{,}}}
\newcommand{\commaNoStar}{\mathchoice
	{,\sbox0{$\mathsurround=0pt\textstyle,$}\kern-3.6\ht0\rlap{\color{Bc}\smash{\rule[-0.999\dp0]{1.4\ht0}{1.1\ht0+1.1\dp0}}}\kern2.021\ht0\llap{\BeginAccSupp{ActualText=}\textup{,}\EndAccSupp{}}\hspace{1.579\ht0}\mspace{3mu plus 0mu minus 3mu}}
	{,\sbox0{$\mathsurround=0pt\textstyle,$}\kern-3.6\ht0\rlap{\color{Bc}\smash{\rule[-0.999\dp0]{1.4\ht0}{1.1\ht0+1.1\dp0}}}\kern2.021\ht0\llap{\BeginAccSupp{ActualText=}\textup{,}\EndAccSupp{}}\hspace{1.579\ht0}\mspace{3mu plus 0mu minus 3mu}}
	{\textup{,}\mspace{2mu plus 0mu minus 2mu}}
	{\textup{,}\mspace{1mu plus 0mu minus 1mu}}}
\renewcommand{\;}{\ecomma}
\newcommand{\ecomma}{\mathchoice
	{,\sbox0{$\mathsurround=0pt\textstyle,$}\kern-3.6\ht0\rlap{\color{Bc}\smash{\rule[-0.999\dp0]{1.4\ht0}{1.1\ht0+1.1\dp0}}}\kern2.021\ht0\llap{\BeginAccSupp{ActualText=},\EndAccSupp{}}}
	{,\sbox0{$\mathsurround=0pt\textstyle,$}\kern-3.6\ht0\rlap{\color{Bc}\smash{\rule[-0.999\dp0]{1.4\ht0}{1.1\ht0+1.1\dp0}}}\kern2.021\ht0\llap{\BeginAccSupp{ActualText=},\EndAccSupp{}}}
	{\text{,}}
	{\text{,}}}
\renewcommand{\:}{\eperiod}
\newcommand{\eperiod}{\mathchoice
	{.\sbox0{$\mathsurround=0pt\textstyle.$}\kern-2\ht0\rlap{\color{Bc}\smash{\rule{1.2\ht0}{1.2\ht0}}}\kern2\ht0\llap{\BeginAccSupp{ActualText=}.\EndAccSupp{}}}
	{.\sbox0{$\mathsurround=0pt\textstyle.$}\kern-2\ht0\rlap{\color{Bc}\smash{\rule{1.2\ht0}{1.2\ht0}}}\kern2\ht0\llap{\BeginAccSupp{ActualText=}.\EndAccSupp{}}}
	{\text{.}}
	{\text{.}}}

\renewcommand{\thickspace}{\tmspace+\thickmuskip{.2777em}}
\makeatother
\usepackage{ulem}
\usepackage{amsthm}
\newtheorem{prop}{Proposition}[section]
\newtheorem*{prop*}{Proposition}
\newtheorem{thm}[prop]{Theorem}
\newtheorem*{thm*}{Theorem}
\newtheorem{lem}[prop]{Lemma}
\newtheorem*{lem*}{Lemma}

\newtheorem*{cor*}{Corollary}
\theoremstyle{remark}
\newtheorem{rmk}{Remark}[prop]
\newtheorem*{rmk*}{Remark}

\newtheorem*{obs*}{Observation}

\newtheorem*{defn*}{Definition}
\theoremstyle{definition}
\newcounter{A}
\usepackage{needspace}
\usepackage{enumerate}
\usepackage{array}
\usepackage{booktabs}
\usepackage{tikz}
\usetikzlibrary{cd}
\usepackage{graphicx}
\usepackage{bookmark}
\usepackage{hyperref}\hypersetup{final,hidelinks,bookmarksnumbered=true,bookmarksopen=true}
\usepackage{macross}
\makeatletter
	\newif\ifstartedinmathmode
	\startedinmathmodefalse
\renewcommand{\start@aligned}[2]{%
	\RIfM@\else\nonmatherr@{\begin{\@currenvir}}\fi\savecolumn@\alignedspace@left\if #1t\vtop \else \if#1b \vbox \else \vcenter \fi \fi\bgroup\maxfields@#2\relax\ifnum\maxfields@>\m@ne\multiply\maxfields@\tw@\let\math@cr@@@\math@cr@@@alignedat\alignsep@\z@skip\else\let\math@cr@@@\math@cr@@@aligned\alignsep@\minalignsep\fi\Let@ \chardef\dspbrk@context\@ne\default@tag\spread@equation\global\column@\z@\ialign\bgroup%
	&\column@plus\hfil\strut@$\m@th\startedinmathmodetrue\textstyle{##{}}\startedinmathmodefalse$\tabskip\z@skip%
	&\column@plus$\m@th\startedinmathmodetrue\textstyle{##}\startedinmathmodefalse$\hfil\tabskip\alignsep@\crcr%
	}
\renewcommand{\align@preamble}{%
	&\hfil\strut@\setboxz@h{\@lign$\m@th\startedinmathmodetrue\textstyle{##{}}\startedinmathmodefalse$}%
	\ifmeasuring@\savefieldlength@\fi\set@field\tabskip\z@skip%
	&\setboxz@h{\@lign$\m@th\startedinmathmodetrue\textstyle{##}\rule{-\leftmarginE}{0.0pt}\startedinmathmodefalse$}%
	\ifmeasuring@\savefieldlength@\fi\set@field\hfil\tabskip\alignsep@%
	}
\renewenvironment{flalign*}{%
	\start@align\tw@\st@rredtrue\m@ne\rule{\leftmarginE}{0.0pt}
	}{%
	\endalign
	}

	\let\originalitemize\itemize
	\let\originalenditemize\enditemize
	\newlength{\leftmarginE}
	\newlength{\leftmarginI}
	\newlength{\leftmarginII}
\renewenvironment{itemize}{%
	\NEWleftmargini{10.00002pt}%
	\NEWleftmarginii{9.25003pt}%
	\NEWleftmarginiii{8.49998pt}%
	\NEWleftmarginiv{7.75002pt}%
	\originalitemize%
	\addtolength{\leftmarginE}{\leftmargin}\tsSwap%
	}{%
	\originalenditemize%
	\tsSwapBack\addtolength{\leftmarginE}{-\leftmargin}%
	\NEWleftmargini{12.7778pt}%
	\NEWleftmarginii{12.7778pt}%
	\NEWleftmarginiii{18.61116pt}%
	\NEWleftmarginiv{15.62502pt}%
	}
\def\@enum@{%
	\list{\csname label\@enumctr\endcsname}{\usecounter{\@enumctr}\def\makelabel##1{\hss\llap{##1}}}%
	\addtolength{\leftmarginE}{\leftmargin}\tsSwap%
	}

\renewenvironment{proof}[1][\proofname]{%
	\par\pushQED{\qed}\normalfont \topsep6\p@\@plus6\p@\relax\trivlist\item\relax{\itshape#1\@addpunct{.}}\hspace\labelsep\ignorespaces%
	\allowdisplaybreaks%
	}{%
	\popQED\endtrivlist\@endpefalse%
	}
	\newlength{\fontdimentwofont}
\newenvironment{iws}[1]{%
	\setlength{\fontdimen2\font}{\fontdimentwofont+#1}\ignorespaces%
	}{%
	\setlength{\fontdimen2\font}{\fontdimentwofont}\ignorespacesafterend%
	}
\makeatother
\begin{document}
\title{Spectral estimates of dynamically-defined and amenable operator families}
\author[S.~Beckus]{Siegfried Beckus}
\address{Institut f\"ur Mathematik, Universit\"at Potsdam, Potsdam~14476, Germany}
\email{beckus@uni-potsdam.de}
\author[A.~Takase]{Alberto Takase}
\address{Department of Mathematics, University of California at Irvine, CA~92697, USA}
\email{atakase@uci.edu}
\begin{abstract}
We consider
	kernel operators defined by a~dynamical system.
The Hausdorff distance of spectra is estimated by the Hausdorff distance of subsystems.
We prove that
	the spectrum map is $ \frac{1}{2} $-Hölder continuous provided
	the group action and kernel are Lipschitz continuous and
	the group has strict polynomial growth.
Also, we prove that
	the continuity can be improved resulting in the spectrum map being Lipschitz continuous provided
	the kernel is instead locally-constant.
This complements a~1990 result by J.~Avron; P.H.M.v.~Mouche; B.~Simon establishing that
	one-dimensional discrete quasiperiodic Schrödinger operators with Lipschitz continuous potentials, e.g., the Almost Mathieu Operator, exhibit spectral $ \frac{1}{2} $-Hölder continuity.
Also, this complements a~2019 result by S.~Beckus; J.~Bellissard; H.~Cornean establishing that
	$ d $-dimensional discrete subshift Schrödinger operators with locally-constant potentials, e.g., the Fibonacci Hamiltonian, exhibit spectral Lipschitz continuity.
Our work exposes the connection between the past two results, and
	the group, e.g., the Heisenberg group, needs not be the integer lattice nor abelian.
\end{abstract}
\date{\today}
\thanks{This project was supported by the Deutsche Forschungsgemeinschaft [\href{https://gepris.dfg.de/gepris/projekt/412141125?language=en}{BE 6789/1-1 to S.B.}] and the second author was partially supported by NSF grant \href{https://www.nsf.gov/awardsearch/showAward?AWD_ID=1855541&HistoricalAwards=false}{DMS-1855541} (PI - A.~Gorodetski) and NSF fellowship award \href{https://www.nsf.gov/awardsearch/showAward?AWD_ID=2213277&HistoricalAwards=false}{DMS-2213277}.}
\maketitle
\section{Introduction}
{
\noindent
The \emph{Almost Mathieu Operator} (AMO) is the linear operator
	$ H_{\lambda,\alpha,\om}:\ell^2(\Z)\to\ell^2(\Z) $
	defined by
	\[\ts
	(H_{\lambda,\alpha,\om}\psi)(n)=\psi(n+1)+\psi(n-1)+2\lambda\cos\p{2\pi(n\alpha+\om)}\psi\p{n}
	\]
	for every $ \psi\in \ell^2(\Z)\, n\in\Z \:$
	Here $ \lambda \;$ $ \alpha\;$ $ \om $ are parameters in $ \R $;
	without loss of generality,
	\[\ts \om\in \T\defeq \R/\Z=[0,1]/{\sim} \:\]
Observe the spectrum is a nonempty compact subset of $ \R $ since the AMO is both bounded and self-adjoint.
The AMO has its origins in solid-state physics and the study of electrons;
	the name was introduced by B.~Simon in 1982 \cite{MP1982_SimonRP}, but
	the operator dates back to 1955 when P.G.~Harper \cite{P1955_Harper1of2,P1955_Harper2of2} under the tutelage of R.E.~Peierls \mbox{published} one of the first descriptions of the spectrum with $ \lambda=1 $ and $ \om=0 $ by utilizing a tight-binding approximation: electrons in a crystal and in a magnetic field have nondiscrete nonevenly-spaced broadened energy values.
Physicists sought to better understand the spectrum;
	M.~Azbel in 1964 \cite{P1964_Azbel} conjectured and D.~Hofstadter in 1976 \cite{P1976_Hofstadter} computationally supported (see Figure~\ref{fig:1}) that the spectrum has the characteristic of being either band-like for rational $ \alpha $ or fractal-like for irrational~$ \alpha \:$
Decades later A.~Avila and S.~Jitomirskaya in 2009 \cite{MP2009_AvilaJitomirskayaAMO} made the final step towards the complete solution of the Ten Martini Problem which sought to confirm the conjectured topological structure of the spectrum: if $ \alpha=\frac*{p}{q} $ is rational, then the spectrum is a disjoint union of at most $ q $-many compact intervals; if $ \alpha $ is irrational, then the spectrum is a Cantor set, i.e., a~nonempty compact subset of $ \R $ not having isolated points nor interior points.
\begin{figure}[t]
\centering
\includegraphics[height=0.334\textheight]{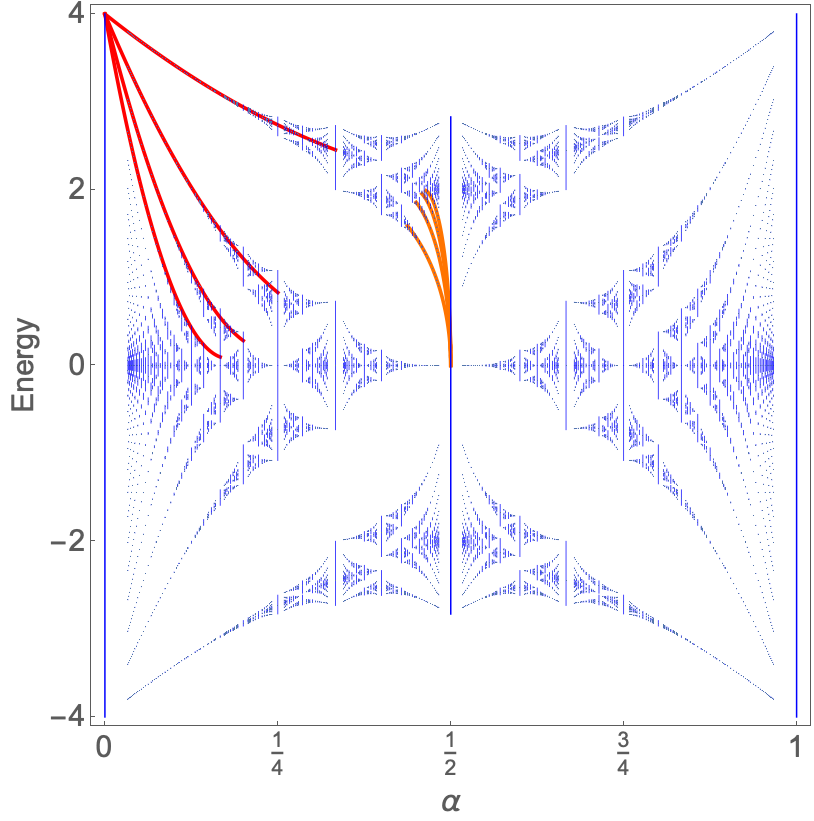}
\caption{\cite{P2016_HatsKatsTach}: Hofstadter's butterfly}
\label{fig:1}
\end{figure}
}
\subsubsection*{Our motivation: the 1990 result}
J.~Avron; P.H.M.v.~Mouche; B.~Simon \cite{MP1990_AvronMoucheSimon} considered a general operator family (one-dimensional discrete quasiperiodic Schrödinger operators) containing the AMO and established that the Hausdorff distance $ {\dist_{\ssm H}} $ of spectra has a square root behavior with respect to the Hausdorff distance $ d_{\ssm H} $ of subsystems or, equivalently, the spectrum map is $ \frac{1}{2} $-Hölder continuous.
Define the one-dimensional torus $ \T\defeq \R/\Z=[0,1]/{\sim} \:$
Let the phase $ \om\in\T $ and let the frequency $ \alpha\in\R \:$
Endow $ \R\times \T $ with the metric
\[\ts
d:(\R\times \T)\times (\R\times \T)\to[0,+\infty):(\b{\begin{subarray}{1}
\alpha_{1}\\
\om_{1}
\end{subarray}},\b{\begin{subarray}{1}
\alpha_{2}\\
\om_{2}
\end{subarray}})\mapsto \max\s{\v{\alpha_1-\alpha_2},\v{\om_1-\om_2}}\:
\]
The (quasiperiodic) dynamical system is $ \ab{\R\times\T,(\Z,\tau)} \;$ i.e., $ \Z $ acts on $ \R\times\T $ with the group action
\[\ts
\tau:\Z\times (\R\times \T)\to \R\times \T: (-n,\b{\begin{subarray}{1}
\alpha\\
\om
\end{subarray}})\mapsto \b{\begin{subarray}{1}
\alpha\\
\bp{T}{\alpha}{-n}\om
\end{subarray}}\;
\]
where $ T_\alpha:\T\to \T:\om\mapsto -\alpha+\om \:$
Let the potential $ v:\T\to\R $ be continuous.
The kernel for the operator family (one-dimensional discrete Schrödinger operators with potential $ v $) is
\[\ts
k_v:\Z\times (\R\times \T)\to\C:(m,\b{\begin{subarray}{1}
\alpha\\
\om
\end{subarray}})\mapsto\begin{cases}
1&\iph m\in\s{-1,1}\\
v(\om)&\iph m=0\\
0&\elss\:
\end{cases}
\]
The operator family is $ H_{v,\alpha,\om}:\ell^2(\Z)\to\ell^2(\Z) $ defined by
\[\ts
(H_{v,\alpha,\om}\psi)(n)=\sum_{m\in \Z}k_v(-n+m,-n\cdot_\tau\b{\begin{subarray}{1}
\alpha\\
\om
\end{subarray}})\psi(m)=\psi(n+1)+\psi(n-1)+v(\bp{T}{\alpha}{-n}\om)\psi(n)
\]
for every $ \psi\in \ell^2(\Z)\, n\in\Z \:$
Observe if $ v(\om)=2\lambda\cos(2\pi\om) \;$ then we recover the AMO.
Endow the collection $ \Ss $ of all left-invariant nonempty closed subsets of $ \R\times\T $ with the Hausdorff metric
\[\ts
d_{\ssm H}:\Ss\times \Ss\to[0,+\infty]:(X,Y)\mapsto \max\s{\sup_{z\in X}d(z,Y),\sup_{z\in Y}d(z,X)}\;
\]
and endow the collection of all nonempty compact subsets of $ \R $ with the Hausdorff metric $ {\dist_{\ssm H}} \:$
The 1990 result \cite{MP1990_AvronMoucheSimon} is the following theorem:
	if $ v:\T\to \R $ is Lipschitz continuous, then there exist $\de>0$ and $ C>0 $ such that for each $ (\alpha_{1},\alpha_{2})\in \R^2\cap\s{\text{\(\de\)-close}} \;$
	\[\ts
	\dist_{\ssm H}(\overline{\bigcup_{\om\in\T}\spec(H_{v,\alpha_{ 1},\om})}
	,\overline{\bigcup_{\om\in\T}\spec(H_{v,\alpha_{ 2},\om})})
	\le C d_{\ssm H}(\s{\alpha_{ 1}}\times \T,\s{\alpha_{ 2}}\times\T)^{\frac{1}{2}}=C\v{\alpha_{ 1}-\alpha_{ 2}}^{\frac{1}{2}}\:
	\]
\subsubsection*{Our motivation: a couple remarks, a sharpness remark, and a list of example applications to come}
We remark that $ \ab{\Z,\#} $ is a second-countable locally-compact Hausdorff group endowed with a left Haar measure (the counting measure).
Also, $ \ab{\R\times\T,(\Z,\tau)} $ is a metrizable space endowed with a left continuous action.
The square root behavior is sharp due to $ \alpha=\frac{1}{2} $ in the middle of Figure~\ref{fig:1}; see \cite{P1990_BellissardRammal,P2016_HatsKatsTach}.
Indeed, for fixed $ n $ such as $ n=0,1,2,3 $ in Figure~\ref{fig:1},
\[\ts
\alpha=\frac{p}{q}\to \frac{1}{2}
\lRightarrow
\sup \spec_{n,0}\p{H_{1,\alpha,0}^{\textsc{amo}}}\asymp\pm\sqrt{16\pi n\v{\alpha-\frac*{1}{2}}-16\pi^2n^2\v{\alpha-\frac*{1}{2}}^2}\;
\]
where $ \spec_{n,E} $ is the $ n $-th band nearby energy $E\;$ but the behavior can be linear such as at $ \alpha=0^+ $
\[\ts
\alpha=\frac{p}{q}\to 0^+
\lRightarrow
\sup \spec_{n,4}\p{H_{1,\alpha,0}^{\textsc{amo}}}\asymp 4-2\pi(2n+1)\alpha+\frac{1}{2}\pi^2(2n^2+2n+1)\alpha^2\:
\]
The 2019 result \cite{MP2019_BeckusBellissardCorneanGCont} and the 2018 result \cite{MP2018_BeckusBellissardDeNittisGCT1of2} motivated us to consider a more general operator family $ \s{A_x}_{x\in Z} $ and establish that the Hausdorff distance $ \dist_{\ssm H} $ of spectra has a quantitative ($ \frac{1}{2} $-Hölder or Lipschitz) behavior with respect to the Hausdorff distance $ d_{\ssm H} $ of subsystems.
We will demonstrate that our unified theory applies for the following special cases:
\begin{itemize}
\item
one-dimensional discrete quasiperiodic Schrödinger operators (known \cite{MP1990_AvronMoucheSimon}).
\item
\(d\)-dimensional discrete subshift Schrödinger operators (known \cite{MP2019_BeckusBellissardCorneanGCont}).
\item
one-dimensional discrete skew-shift Schrödinger operators (new).
\item
one-dimensional discrete limit-periodic Schrödinger operators (new).
\end{itemize}
We remark that $ \Z^d $ is special.
One could consider  lattices in a locally compact (Lie) group; see \cite{MP2018_BeckusBellissardDeNittisGCT1of2}.
\subsubsection*{Our work: amenability and continuity}
Let $ \ab{G,\lambda} $ be
	a second-countable locally-compact Hausdorff group endowed with
	a left Haar measure;
	let $ \ell $ be a left-invariant proper $ \Tc $-generating metric on $ G $ and
	define, for $ g\in G \;$ $ \v{g}\defeq \ell(g,e) $ and
	define, for $ r\ge 0 \;$ $ B(r)\defeq \s{g\in G:\v{g}<r} \:$
Let $ \ab{Z,(G,\alpha)} $ be
	a metrizable space endowed with
	a left continuous action;
	let $ d $ be a $ \Tc $-generating metric on $ Z \:$
Let $ k:G\times Z\to\C \;$
	where $ k $ is measurable and
	$\sup_{x\in Z}\int \vv{k(h,\alpha(\inv{(\cdot)},x))}_{\infty} dh<+\infty\:$
Define, for $ x\in Z \;$
	${A}_x:L^2(G)\to L^2(G)$ by
\[\ts
(A_x\psi)(g)=\int k(\inv{g}h,\inv{g}\cdot_\alpha x)\psi(h)dh
\]
for every $ \psi\in L^2(G)\, g\in G \:$
Let $ \Ss $ be the collection of all subsystems of $ \ab{Z,(G,\alpha)} \:$
Define, for $ X\in\Ss \;$
	$ \spec(A_X)\defeq \overline{\bigcup_{x\in X}\spec(A_x)} \:$
The dynamically-defined operators $ A_x $ and their spectra are the main objects of study in this paper.
Our attention is on the regularity of the spectrum map $ X\mapsto \spec\p{A_X} \:$
Using standard techniques such as working with approximate eigenvectors and utilizing strong continuity, we immediately obtain that if $ \ab{G,\lambda} $ is unimodular and $ A_x $ is normal for every $ x\in Z $ and $ k $ satisfies both a continuity condition and a decay condition, then $\spec\p{A_x}=\spec\p{A_{X\p{x}}}\;$ where $ X\p{x}\defeq \overline{\s{\inv{g}x:g\in G}} $ (orbit-closure).
We forgo describing the conditions on $ k $ in detail, but they can be found within \textsc{section~\ref{sec:A1}}.
It suffices to say that stricter conditions on $ k $ are sufficient to prove the following proposition regarding the spectrum map.
\begin{prop}
Assume $ \ab{G,\lambda} $ is unimodular, $ A_x $ is normal for every $ x\in Z \;$ and the following.
\begin{enumerate}[(i)]
\item
$ \alpha $ satisfies a uniform continuity condition:
	for each $ \ep>0 $ and for each nonempty compact subset $ K $ of $ G \;$ there exists $ \de>0 $ such that
	$ d(y,x)<\de \lRightarrow \sup_{g\in K}d(\inv{g}y,\inv{g}x)<\ep \:$
\item
$ k $ satisfies a uniform continuity condition:
	for each $ \ep>0 $ and for each nonempty compact subset $ F $ of $ G \;$ there exists $ \de>0 $ such that
	$ d(y,x)<\de \lRightarrow \int\v{k(h,y)-k(h,x)}\1_F\p{h}dh<\ep\:$
\item
$ k $ satisfies a uniform decay condition:
	for each $ \ep>0 \;$ there exists a nonempty compact subset $ F $ of $ G $ such that
	$ \sup_{x\in Z}\int\vv{k(h,\inv{(\cdot)}x)\1_{G\setminus F}(h)}_{\infty}dh<\ep \:$
\end{enumerate}
Fix $ X\in\Ss \:$
Then
\[\ts \lim_{Y\xrightarrow[]{\smash{\Ss}}X}\sup_{E\in\spec(A_X)}\dist(E,\spec(A_Y))=0\: \]
\end{prop}
\noindent
More can be obtained provided the group $ \ab{G,\lambda} $ is amenable; see \textsc{section~\ref{subsubsec:a}} for two definitions.
This is the context of the following theorem which establishes that the spectrum map is continuous.
\begin{thm}\label{C}
Assume $ \ab{G,\lambda} $ is unimodular, $ A_x $ is normal for every $ x\in Z \;$ and the following.
\begin{enumerate}[(i)]
\item\label{Ca1}
$ \alpha $ satisfies a uniform continuity condition:
	for each $ \ep>0 $ and for each nonempty compact subset $ K $ of $ G \;$ there exists $ \de>0 $ such that
	$ d(y,x)<\de \lRightarrow \sup_{g\in K}d(\inv{g}y,\inv{g}x)<\ep \:$
\item\label{Ck1}
$ k $ satisfies a uniform continuity condition:
	for each $ \ep>0 $ and for each nonempty compact subset $ F $ of $ G \;$ there exists $ \de>0 $ such that
	$ d(y,x)<\de \lRightarrow {\int\v{k(h,y)-k(h,x)}\1_F\p{h}dh}<\ep\:$
\item\label{Ck2}
$ k $ satisfies a uniform decay condition:
	for each $ \ep>0 \;$ there exists a nonempty compact subset $ F $ of $ G $ such that
	$ \sup_{x\in Z}\int\vv{k(h,\inv{(\cdot)}x)\1_{G\setminus F}(h)}_{\infty}dh<\ep \:$
\item\label{Cg1}
$ \ab{G,\lambda} $ is amenable.
\end{enumerate}
Then
\[\ts
\lim_{\de\to 0^+}\sup_{X,Y\in\Ss,d_{\ssm{H}}(X,Y)< \de}\dist_{\ssm{H}}(\spec(A_X),\spec(A_Y))=0\:
\]
\end{thm}
\begin{rmk}
\leavevmode%
\begin{enumerate}[(a)]
\item
In the AMO-like case where $ \ab{G,\lambda}=\ab{\Z,\#} $ and $ \ab{Z,(G,\alpha)}=\ab{\R\times \T,(\Z,\tau)} $ and $ A_{k_v,x}=H_{v,x} \;$ the spectrum map is continuous due to Theorem~\ref{C}.
Indeed, the group $ \ab{\Z,\#} $ is both unimodular and amenable and the potential $ v $ is continuous.
We remark that this corollary---for this specific case---is known; publications of this prior result in the 1980s can be found on page 258 within \cite{MP1982_ElliottAMO} and on page 382 within \cite{MP1983_AvronSimonAMO} and on page 2205 within \cite{MP1985_AvronSimonGCont}.
\item
Theorem~\ref{C} intersects a 2018 result by S.~Beckus; J.~Bellissard; G.~De~Nittis \cite{MP2018_BeckusBellissardDeNittisGCT1of2} establishing that the spectrum map is continuous.
Indeed, the étale groupoid within \cite{MP2018_BeckusBellissardDeNittisGCT1of2} is assumed to satisfy an amenability condition.
We remark that this prior result implies Theorem~\ref{C}; our proof is different and significant since Theorem~\ref{H}---$ \gamma $-Hölder continuity---is obtainable.
\item
If the reader is comfortable with the phase space $ Z $ being compact, then \textit{(\ref{Ca1})} immediately holds.
If the reader is comfortable with the potential $ v $ being discontinuous in the AMO-like case, then \textit{(\ref{Ck1})} needs not hold thus $ x\mapsto H_{v,x}\psi $ needs not be continuous.
We remark that models of quasicrystals have potentials defined on $ \T $ being indicator functions $ \1_{[1-\alpha,1)} $ in the AMO-like case thus discontinuous in the euclidean topology, and the spectrum map is discontinuous at rational frequencies $ \alpha \;$ but in 1991 \cite{MP1991_BellissardIochumTestardGCont} it was established (and plotted) that the spectrum map is continuous at irrational frequencies $ \alpha \;$ where---in a different perspective---the potential defined on $ \s{0,1}^\Z $ is `evaluation at zero' thus continuous in the product topology; see \textsc{section~\ref{subsec:mdsSo}}.
If~the reader is comfortable with the kernel $ k $ having compact range, e.g., $ k_v $ has range $ \s{-1,0,1} $ in the AMO-like case or the self-adjoint Jacobi operator case \cite{MP2000_TeschlReview}, then \textit{(\ref{Ck2})} immediately holds.
If the reader is comfortable with the group $ G $ being $ \Z^d $ or $ \R^d $ or compact or abelian or the Heisenberg group or the lamplighter group, then $ \ab{G,\lambda} $ is both unimodular and amenable \textit{(\ref{Cg1})}.
\end{enumerate}
\end{rmk}
\noindent
More can be obtained provided stricter conditions replace \textit{(\ref{Ca1}--\ref{Cg1})}.
This is the context of the following main theorem which establishes that the spectrum map is $ \gamma $-Hölder continuous ($ \gamma=\frac{1}{2}\text{ or }1 $).
\begin{thm}\label{H}
Assume $ \ab{G,\lambda} $ is unimodular, $ A_x $ is normal for every $ x\in Z \;$ and the following.
\begin{enumerate}[(i)]
\item\label{Ha1}
$ \alpha $ satisfies a Lipschitz continuity condition:
	there exists $ c_\alpha>0 $ such that for each $ g\in G \;$
	\nobreak\hfil\allowbreak\hfilneg
	$ d(\inv{g}x,\inv{g}y)\le (c_\alpha\v{g}+1)d(x,y)\: $
\item\label{Hk1}
$ k $ satisfies one of two quantitative continuity conditions considered; see part \textit{(\ref{Hk1a})} and \textit{(\ref{Hk1b})} below.
\item\label{Hk2}
$ k $ satisfies a linear decay condition:
	there exists $ c_s>0 $ such that for each $ r\ge 0 \;$
	\nobreak\hfil\allowbreak\hfilneg
	$ \sup_{x\in Z}\int \vv{k(h,\inv{(\cdot)}x)\min\s{\v{h},r}}_\infty dh\le c_s \: $
\item\label{Hg1}
$ \ab{(G,\ell),\lambda} $ is strictly-polynomial-growing:
	there exist $ b\ge 1$ and $ c_1\ge c_0>0 $ such that for each $ r\ge 1 \;$
	\nobreak\hfil\allowbreak\hfilneg
	$c_0r^b\le \lambda\p{B(r)}\le c_1r^b \:$
\end{enumerate}
Then
\begin{enumerate}[(a)]
\item\label{Hk1a}
If $ k $ satisfies a Lipschitz continuity condition:
	there exists $ c_k\in L^1_+(G) $ such that for $ \lambda $-a.e.\ $ h\in G \;$
	$\v{k(h,x)-k(h,y)}\le c_{k}\p{h}d(x,y)\;$ then for all $ X,Y\in\Ss \;$
\[\ts
d_{\ssm H}(X,Y)\le{\min\s{\de,1}}
\lRightarrow
\dist_{\ssm H}\p{\spec\p{A_X},\spec\p{A_Y}}\le C d_{\ssm H}(X,Y)^{\frac*{1}{2}}\;
\]
where $ \de\defeq \frac*{c_s}{\vv{c_k}_1c_{\alpha}}(\frac*{(2+b)^{2+b}c_1}{2b^bc_0})^{\frac*{1}{2}} $ and $ C\defeq 2(\vv{c_k}_1c_{\alpha}c_s(\frac*{(2+b)^{2+b}c_1}{2b^bc_0})^{\frac*{1}{2}})^{\frac*{1}{2}}+\vv{c_k}_1 \:$
\item\label{Hk1b}
If $ k $ satisfies a locally-constant condition:
	there exists $ c_k\in L^\infty_+(G) $ such that for $ \lambda $-a.e.\ $ h\in G \;$
	$ d(x,y)<\frac*{1}{c_{k}\p{h}}\lRightarrow k(h,x)=k(h,y)\; $ then for all $ X,Y\in\Ss \;$
\[\ts
d_{\ssm H}(X,Y)\le \de\lRightarrow
\dist_{\ssm H}\p{\spec\p{A_X},\spec\p{A_Y}}\le C d_{\ssm H}(X,Y)\;
\]
where $ \de\defeq \frac*{1}{(\vv{c_k}_\infty+1)(c_\alpha+1)}$ and $ C\defeq (\vv{c_k}_\infty+1)(c_\alpha+1)c_s(\frac*{(2+b)^{2+b}c_1}{2b^bc_0})^{\frac*{1}{2}} \:$
\end{enumerate}
\end{thm}
\begin{rmk}
\leavevmode%
\begin{enumerate}[(a)]
\item
Theorem~\ref{H}(\ref{Hk1a}) generalizes the 1990 result by J.~Avron; P.H.M.v.~Mouche; B.~Simon \cite{MP1990_AvronMoucheSimon}, and
Theorem~\ref{H}(\ref{Hk1b}) generalizes the 2019 result by S.~Beckus; J.~Bellissard; H.~Cornean \cite{MP2019_BeckusBellissardCorneanGCont}.
Our work exposes the connection between the past two results and reveals that quantitative continuity of the spectrum map follows from quite general operator properties.
\item
Lipschitz continuity of the action \textit{(\ref{Ha1})} and strict polynomial growth of the group \textit{(\ref{Hg1})} are specific conditions covering a large class of models.
We remark that our proof allows one to weaken or strengthen both conditions.
For example, we will demonstrate that our unified theory applies in \textsc{section~\ref{subsec:odlSo}} after condition \textit{(\ref{Ha1})} is weakened so that the group action is instead assumed to be exponential Lipschitz leading to a log-Hölder continuity result.
Also, one can generalize Lemma~\ref{lem:po(5)} (where the polynomial growth-estimates of the group, e.g., the Heisenberg group, are utilized) so that yet another quantitative continuity result can be obtained for groups with (sub)exponential growth-estimates.
Specifically, Lemma~\ref{lem:po(5)} is a machine that has as an input the polynomial growth-estimates and that has as an output a suitable family of almost invariant vectors to be utilized later as cutoff functions.
The generalization reduces to replacing growth-estimates of the group with estimates from the suitable family of almost invariant vectors.
We remark that the lamplighter group has exponential growth.
\item
Theorem~\ref{H}(\ref{Hk1a}) is sharp due to the 1990 works \cite{P1990_BellissardRammal,MP1990_HelfferSjostrandAMO}; see also \cite{MP1994_BellissardGCont}.
We remark that
	the regularity can be improved
	to be Lipschitz
	provided the kernel is instead locally-constant; see Theorem~\ref{H}(\ref{Hk1b}).
Also,
	in 1994 \cite{MP1994_BellissardGCont} it was established that
	the regularity at the boundary of a spectral gap can be improved
	to be Lipschitz for a general operator family containing the AMO
	provided the gap width is locally-positive.
Also,
	in 1998 \cite{MP1998_JitomirskayaLastAMO} it was established that
	the regularity can be improved
	to be Lipschitz-log for the AMO
	provided the coupling $ \lambda $ is in the open set $ (0,\frac{2}{29})\cup(\frac{29}{2},+\infty) $; see also \cite{MP2002_JitomirskayaKrasAMO,MP2020_ZhaoAMO}.
Also,
	in 2014 \cite{MP2014_JitomirskayaMaviAMO} it was established that
	the regularity can be generalized
	to be $ \frac{\gamma}{\gamma+1} $-Hölder for one-dimensional discrete quasiperiodic Schrödinger operators with $ \gamma $-Hölder continuous potentials;
	our proofs utilize almost invariant vectors as cutoff functions which is inspired by the 1990 works \cite{MP1990_AvronMoucheSimon,MP1990_ChoiElliottYui} which itself inspired the generalization within \cite{MP2014_JitomirskayaMaviAMO}, and we remark that
	our regularity can too be generalized
	to be $ \frac{\gamma}{\gamma+1} $-Hölder in the setting for Theorem~\ref{H}(\ref{Hk1a})
	provided the kernel is instead $ \gamma $-Hölder continuous.
\end{enumerate}
\end{rmk}
\subsubsection*{Our work: final remarks}
To explain why the continuous behavior of the spectra has been studied for various models, there are cases where one has a sequence of (periodic) approximations and thus a sequence of convergent spectra.
This has been useful to analyze the spectral theory of operators.
For example, rational approximations provide a deeper understanding of the Almost Mathieu Operator and related quasiperiodic models.
Similarly, periodic approximations provide a deeper understanding of the Fibonacci Hamiltonian and related Sturmian models.
For more papers about `almost periodic' operators, we recommend the 1990 paper by A.~Süt\H{o} \cite{MP1990_SutoFH} (also \cite{MP1992_BellissardGLT}) and the references therein.
Also, we recommend the 2016 paper by D.~Damanik; A.~Gorodetski; W.~Yessen \cite{MP2016_DamanikGorodetskiYessFH}.
Approximations are a practical tool to compute the spectrum.
Indeed, D.~Hofstadter in 1976 \cite{P1976_Hofstadter} plotted one such spectral approximation for the AMO (see Figure~\ref{fig:1}), but what is more important is that it can be utilized to derive spectral properties of the limit-operator $ A_x $ from the periodic approximants $ A_{y_j} \:$
\par
Indeed, physicists sought to better understand electrons in a crystal and in a magnetic field;
	E.J.~Austin and M.~Wilkinson in 1990 \cite{P1990_AustinWilkinson} computationally studied a general operator family containing the AMO and plotted a Hofstadter-like butterfly.
Also, physicists sought to better understand the fractal dimension of $ \spec\p{\bp{H}{1,\alpha,0}{\textsc{amo}}} $;
	M.~Kohmoto and C.~Tang in 1986 \cite{P1986_KohmotoTang} computationally studied $ \alpha=(\sqrt{5}-1)/2\, (\sqrt{2}-1)\, (\pi-3)\sqrt{3} $
	and conjectured that the dimension is~$ \frac{1}{2} $ for irrational~$ \alpha $ (the dimension is~$ 1 $ for rational~$ \alpha $); see also \cite{P1987_BellStinchcombe,P1991_GeiselKetzmerickPetschel}.
Years later 
	Y.~Last in 1994 \cite{MP1994_LastAMO} proved that the Hausdorff dimension does not exceed~$ \frac{1}{2} $ for topologically-generic irrational~$ \alpha^{({\sss \Tc})} $---the proof utilizes~$ \frac{1}{2} $-Hölder regularity of the spectrum map.
But,
	E.J.~Austin and M.~Wilkinson in 1994 \cite{P1994_AustinWilkinson} proposed that the box-counting dimension is~nonconstant for irrational~$ \alpha $ and approaches zero for specific~$ \alpha_n \;$ and
	T.~Geisel; R.~Ketzmerick; K.~Kruse; F.~Steinbach in 1998 \cite{P1998_GeiselKetzmerickKruseSteinbach} proposed that, given an assumption, the Hausdorff dimension does not exceed~$ \frac{1}{2} $ for all irrational~$ \alpha \:$
A~decade later
	S.~Jitomirskaya and I.~Krasovsky in 2019 \cite{MP2019_JitomirskayaKrasAMO} made the final step towards the conjectured fractal dimension of the spectrum: the Hausdorff dimension does not exceed~$ \frac{1}{2} $ for \emph{all} irrational~$ \alpha $---the proof utilizes Lipschitz-log regularity of the spectrum map.
We remark that
	Y.~Last and M.~Shamis in 2016 \cite{MP2016_LastShamisAMO} proved that the Hausdorff dimension is~$0$ for topologically-generic irrational~$ \alpha^{({\sss \Tc'})} $---the proof utilizes~$ \frac{1}{2} $-Hölder regularity of the spectrum map---and
	S.~Jitomirskaya and S.~Zhang in 2022 \cite{MP2022_JitomirskayaZhangAMO} proved that the box-counting dimension is~$1$ for the same topologically-generic irrational~$ \alpha^{({\sss \Tc'})} $ and $ \Tc'\subseteq \Tc \:$
Also,
	B.~Helffer; Q.~Liu; Y.~Qu; Q.~Zhou in 2019 \cite{MP2019_HelfferLiuQuZhouAMO} proved that the collection $ \Tc'' $ of all irrational $ \alpha $ such that the Hausdorff dimension is~positive (thus $ \Tc''\cap\Tc'=\es $) is~topologically-dense and itself has positive Hausdorff dimension.
\par
Within \cite{MP2020_AvniBreuerSimonLGO},
	periodic Jacobi matrices on trees are studied, and
	our theorem does not apply since
	we consider groups but they consider trees.
We remark that the \cite{MP2020_AvniBreuerSimonLGO} authors themselves remark that
	what makes their objects fascinating is that the fundamental group is nonabelian.
Within \cite{MP2022_BrunoCalziLGO},
	Schrödinger operators on Lie groups are studied, and
	our theorem does not apply since
	we consider dynamically-defined operators but they do not have a dynamical system.
	Also, their group needs not be unimodular thus needs not have polynomial growth but their group satisfies a weak polynomial growth condition.
We remark that the \cite{MP2022_BrunoCalziLGO} authors themselves remark that
	what makes their theorem fascinating is that their work applies to the Heisenberg group (a nonabelian group).
Within \cite{MP2007_LenzPeyerVeGAmen},
	random Schrödinger operators on manifolds are studied.
We remark that the \cite{MP2007_LenzPeyerVeGAmen} authors themselves remark that
	their work continues the \cite{MP1985_BellissardLimaTestard} approach (see also \cite{MP1986_Bellissard}), and
	what makes their theorem fascinating is that they use techniques from Connes' noncommutative integration theory and von~Neumann algebras.
For more related papers, see \cite{MR2962855,MR4366010}.
\par
Within \cite{MP1985_BellissardLimaTestard},
	our operator families are studied, and
	we consider spectral continuity results but they consider the integrated density of states and spectral gap-labeling results.
We remark that the \cite{MP1985_BellissardLimaTestard} authors themselves remark that
	their work generalizes results by both R.~Johnson; J.~Moser (1982) and F.~Delyon; B.~Souillard (1983), and
	what makes their theorem fascinating is that they can deal with an arbitrary dimension both in the discrete and in the continuous case.
\par
In \textsc{section~\ref{sec:ea}} we demonstrate example applications.
In \textsc{section~\ref{sec:p}} we provide the preliminaries.
In \textsc{section~\ref{sec:po}} we prove the main results.
In \textsc{section~\ref{sec:A1}} we prove a standard lemma.
\clearpage
\section{Example Applications}\label{sec:ea}
\subsection{One-dimensional Discrete Quasiperiodic Schrödinger Operators}\label{subsec:odqSo}
We demonstrate that our unified theory applies for this titular setting.

Let $ \ab{\Z,\#} $ be the one-dimensional lattice endowed with the counting measure.
Let $ \ell $ be the metric on $ \Z $ defined by
\[\ts
\ell(n,m)=\v{-n+m}
\]
for every $ n,m\in\Z \:$
Define the phase space $ Z\defeq \R\times\T \:$
Let $ d $ be the metric on $ Z $ defined by
\[\ts
d(x,y)=\max\s{\v{x_1-y_1},\v{x_2-y_2}}
\]
for every $ x,y\in Z \:$
Let $ \ab{Z,(\Z,\tau)} $ be the (quasiperiodic) dynamical system defined by
	\[\ts
	-n\cdot_{\tau} x=
	\b*{\begin{subarray}{1}
	x_1\\
	nx_1+x_2
	\end{subarray}}
	\]
	for every $ n\in\Z\, x\in Z \:$
Define, for $ v\in C(\T) \;$ the kernel $ k_v=k:\Z\times Z\to\C $ by
	\[\ts
	k(m,x)=\begin{cases}
	1&\iph m\in\s{-1,1}\\
	v(x_2)&\iph m=0\\
	0&\elss
	\end{cases}
	\]
	for every $ m\in\Z\, x\in Z \:$
Define the operator family $ \s{A_{v,x}}_{v\in C(\T),x\in Z}\subseteq \Ls\p{\ell^2(\Z)} $ by
	\[\ts
	(A_{x}\psi)(n)=\sum_{m\in\Z}k(-n+m,-n\cdot_\tau x)\psi(m)=\psi(n+1)+\psi(n-1)+v(nx_1+x_2)\psi(n) \]
	for every $ v\in C(\T)\, x\in Z\, \psi\in \ell^2(\Z)\, n\in\Z \:$
Observe if $ v\p{\om}=2\lam \cos\p{2\pi \om} \;$ then we recover the Almost Mathieu Operator.
\begin{thm}[$ \frac{1}{2} $-Hölder regularity]\label{thm:odqSo1}
Assume the potential $ v:\T\to\R $ is Lipschitz continuous, i.e., there exists $ \rrm c_v=\rrm c>0 $ such that for all $ \om_1,\om_2\in\T \;$ $ \v{v(\om_1)-v(\om_2)}\le \rrm c\mspace{1mu minus 1mu}\v{\om_1-\om_2} \:$
Then
\begin{enumerate}[(a)]
\item\label{odqSo1a}
$ \tau $ satisfies a Lipschitz continuity condition:
	for each $ n\in\Z \;$
	\[\ts
	d(-n\cdot_\tau x,-n\cdot_\tau y)\le \p{\v{n}+1}d(x,y)\:
	\]
\item\label{odqSo1b}
$ k $ satisfies a Lipschitz continuity condition:
	for each $ m\in \Z \;$
	\[\ts \v{k(m,x)-k(m,y)}\le \rrm c\mspace{1mu minus 1mu}\de_0(m)d(x,y)\: \]
\item\label{odqSo1c}
$ k $ satisfies a linear decay condition:
	for each $ r\ge 0 \;$
	\[\ts
	\sup_{x\in Z}\sum_{m\in\Z}\vv{k(m,-(\cdot)\cdot_\tau x)\min\s{\v{m},r}}_\infty\le 2\:
	\]
\item\label{odqSo1d}
$ \ab{(\Z,\ell),\#} $ is strictly-polynomial-growing:
	for each $ r\ge 1 \;$
	\[ r\le \# B(r)\le 3r \:\]
\item\label{odqSo1e}
For all $ X,Y\in\Ss \;$
\[\ts
d_{\ssm H}(X,Y)\le{\min\s{\de,1}}
\lRightarrow
\dist_{\ssm H}\p{\spec\p{A_X},\spec\p{A_Y}}\le C d_{\ssm H}(X,Y)^{\frac*{1}{2}}\;
\]
where $ \de\defeq \frac*{2}{\rrm c}(\frac*{3^4}{2})^{\frac*{1}{2}} $ and $ C\defeq 2(\rrm c\mspace{1mu minus 1mu}\frac{2}{1}(\frac*{3^4}{2})^{\frac*{1}{2}})^{\frac*{1}{2}}+\rrm c \:$
\end{enumerate}
\end{thm}
\noindent
Theorem~\ref{thm:odqSo1}(\ref{odqSo1a}, \ref{odqSo1b}, \ref{odqSo1c}, \ref{odqSo1d}) are easy to see and come in handy since Theorem~\ref{thm:odqSo1}(\ref{odqSo1e}) follows from a combination of (\ref{odqSo1a}--\ref{odqSo1d}) and Theorem~\ref{H}(\ref{Hk1a}).
We remark that a 1990 result by both \cite{MP1990_ChoiElliottYui} and \cite{MP1990_AvronMoucheSimon}
	established that the spectrum map is Hölder continuous in this specific setting;
	the regularity was established to be $ \frac{1}{3} $-Hölder \cite{MP1990_ChoiElliottYui} but was improved to be $ \frac{1}{2} $-Hölder \cite{MP1990_AvronMoucheSimon}.
Our estimate is different and significant since $ \de=\frac*{2}{\rrm c}(\frac*{3^4}{2})^{\frac*{1}{2}} $ is explicitly written; another $ C $ constant is explicitly written within \cite{MP1990_AvronMoucheSimon}.
\subsection{Multi-dimensional Discrete Subshift Schrödinger Operators}\label{subsec:mdsSo}
We demonstrate that our unified theory applies for this titular setting.

Let $ \ab{\Z^d,\#} $ be the multi-dimensional lattice endowed with the counting measure.
Let $ \ell $ be the metric on $ \Z^d $ defined by
\[\ts
\ell(n,m)=\max\{\v{-n_1+m_1},\ldots,\v{-n_d+m_d}\}
\]
for every $ n,m\in\Z^d \:$
Define the phase space $ Z\defeq \s{0,1}^{\Z^d} \:$
Let $ d $ be the metric on $ Z $ defined by
\[\ts
d(x,y)=\inf\s{\frac{1}{r+1}:r\in[0,+\infty)\, x|_{B(r)}=y|_{B(r)}}
\]
for every $ x,y\in Z \:$
Let $ \ab{Z,(\Z^d,\tau)} $ be the (subshift) dynamical system defined by
	\[\ts
	-n\cdot_{\tau} x=
	(m\mapsto x(n+m))\eqdef T^{-n}x
	\]
	for every $ n\in\Z^d\, x\in Z \:$
Define, for $ v\in C(Z) \;$ the kernel $ k_v=k:\Z^d\times Z\to\C $ by
	\[\ts
	k(m,x)=\begin{cases}
	1&\iph \vv{m}_1=1\\
	v(x)&\iph m=\vec{0}\\
	0&\elss
	\end{cases}
	\]
	for every $ m\in\Z^d\, x\in Z \:$
Define the operator family $ \s{A_{v,x}}_{v\in C(Z),x\in Z}\subseteq \Ls\p{\ell^2(\Z^d)} $ by
	\[\ts
	(A_{x}\psi)(n)=\sum_{m\in\Z^d}k(-n+m,-n\cdot_\tau x)\psi(m)=\p{\sum_{\vv{m}_1=1}\psi(n+m)}+v(T^{-n}x)\psi(n) \]
	for every $ v\in C(Z)\, x\in Z\, \psi\in \ell^2(\Z^d)\, n\in\Z^d \:$
Observe if $ v\p{x}=x\p{\vec{0}} \;$ then we recover the Fibonacci Hamiltonian \cite{MP2016_DamanikGorodetskiYessFH}.
\begin{thm}[Lipschitz regularity]\label{thm:mdsSo1}
Assume the potential $ v:Z\to\R $ is locally-constant, i.e., there exist $ \rrm r_v=\rrm r\ge 0 $ and $ \rrm v_v=\rrm v:\s{0,1}^{B(\rrm r)}\to \R $ such that for each $ x\in Z \;$ $ v(x)=\rrm v(x|_{B(\rrm r)}) \:$
Then
\begin{enumerate}[(a)]
\item\label{mdsSo1a}
$ \tau $ satisfies a Lipschitz continuity condition:
	for each $ n\in\Z^d \;$
	\[\ts
	d(-n\cdot_\tau x,-n\cdot_\tau y)\le \p{\v{n}+1}d(x,y)\:
	\]
\item\label{mdsSo1b}
$ k $ satisfies a locally-constant condition:
	for each $ m\in \Z^d \;$
	\[\ts d(x,y)<\frac{1}{(\rrm r+1)\de_0(m)}\lRightarrow
	k(m,x)=k(m,y)\: \]
\item\label{mdsSo1c}
$ k $ satisfies a linear decay condition:
	for each $ r\ge 0 \;$
	\[\ts
	\sup_{x\in Z}\sum_{m\in\Z^d}\vv{k(m,-(\cdot)\cdot_\tau x)\min\s{\v{m},r}}_\infty\le 2d\:
	\]
\item\label{mdsSo1d}
$ \ab{(\Z^d,\ell),\#} $ is strictly-polynomial-growing:
	for each $ r\ge 1 \;$
	\[ r^d\le \# B(r)\le 3^dr^d \:\]
\item\label{mdsSo1e}
For all $ X,Y\in\Ss \;$
\[\ts\hspace*{\leftmargin minus \leftmargin}
d_{\ssm H}(X,Y)\le \de\lRightarrow
\dist_{\ssm H}\p{\spec\p{A_X},\spec\p{A_Y}}\le C d_{\ssm H}(X,Y)\;
\]
where $ \de\defeq \frac*{1}{2(\rrm r+2)}$ and $ C\defeq 2(\rrm r+2)\frac{2d}{1}(\frac*{(2+d)^{2+d}3^d}{2d^d})^{\frac*{1}{2}} \:$
\end{enumerate}
\end{thm}
\noindent
Theorem~\ref{thm:mdsSo1}(\ref{mdsSo1a}, \ref{mdsSo1b}, \ref{mdsSo1c}, \ref{mdsSo1d}) are easy to see and come in handy since Theorem~\ref{thm:mdsSo1}(\ref{mdsSo1e}) follows from a combination of (\ref{mdsSo1a}--\ref{mdsSo1d}) and Theorem~\ref{H}(\ref{Hk1b}).
We remark that a 2019 result by \cite{MP2019_BeckusBellissardCorneanGCont}
	established that the spectrum map is Lipschitz continuous in this specific setting.
Our estimate is different and significant since $ \de=\frac*{1}{2(\rrm r+2)} $ is explicitly written; another $ C $ constant is explicitly written within \cite{MP2019_BeckusBellissardCorneanGCont}.
\subsection{One-dimensional Discrete Skew-shift Schrödinger Operators}\label{subsec:odsSo}
We demonstrate that our unified theory applies for this titular setting provided a lenient condition replaces condition \textit{(\ref{Ha1})} within Theorem~\ref{H}; see Theorem~\ref{thm:odsSo1}(\ref{odsSo1a}).

Let $ \ab{\Z,\#} $ be the one-dimensional lattice endowed with the counting measure.
Let $ \ell $ be the metric on $ \Z $ defined by
\[\ts
\ell(n,m)=\v{-n+m}
\]
for every $ n,m\in\Z \:$
Define the phase space $ Z\defeq \R\times\T\times\T \:$
Let $ d $ be the metric on $ Z $ defined by
\[\ts
d(x,y)=\max\s{\v{x_1-y_1},\v{x_2-y_2},\v{x_3-y_3}}
\]
for every $ x,y\in Z \:$
Let $ \ab{Z,(\Z,\tau)} $ be the (skew-shift) dynamical system defined by
	\[\ts
	-n\cdot_\tau x=
	\b*{\begin{subarray}{1}
	x_1\\
	nx_1+x_2\\
	\frac*{n(n-1)}{2}x_1+nx_2+x_3
	\end{subarray}}
	\]
	for every $ n\in\Z\, x\in Z \:$
Define, for $ v\in C(\T) \;$ the kernel $ k_v=k:\Z\times Z\to\C $ by
	\[\ts
	k(m,x)=\begin{cases}
	1&\iph m\in\s{-1,1}\\
	v(x_3)&\iph m=0\\
	0&\elss
	\end{cases}
	\]
	for every $ m\in\Z\, x\in Z \:$
Define the operator family $ \s{A_{v,x}}_{v\in C(\T),x\in Z}\subseteq \Ls\p{\ell^2(\Z)} $ by
	\[\ts
	(A_{x}\psi)(n)=\sum_{m\in\Z}k(-n+m,-n\cdot_\tau x)\psi(m)=\psi(n+1)+\psi(n-1)+v(\frac{n(n-1)}{2}x_1+nx_2+x_3)\psi(n) \]
	for every $ v\in C(\T)\, x\in Z\, \psi\in \ell^2(\Z)\, n\in\Z \:$
\begin{thm}[$ \frac{1}{3} $-Hölder regularity]\label{thm:odsSo1}
Assume the potential $ v:\T\to\R $ is Lipschitz continuous, i.e., there exists $ \rrm c_v=\rrm c>0 $ such that for all $ \om_1,\om_2\in\T \;$ $ \v{v(\om_1)-v(\om_2)}\le \rrm c\mspace{1mu minus 1mu}\v{\om_1-\om_2} \:$
Then
\begin{enumerate}[(a)]
\item\label{odsSo1a}
$ \tau $ satisfies a power Lipschitz continuity condition:
	for each $ n\in \Z \;$
	\[\ts d(-n\cdot_\tau x,-n\cdot_\tau y)\le (\frac{n(n-1)}{2}+\v{n}+1)d(x,y)\: \]
\item\label{odsSo1b}
$ k $ satisfies a Lipschitz continuity condition:
	for each $ m\in \Z \;$
	\[\ts \v{k(m,x)-k(m,y)}\le \rrm c\mspace{1mu minus 1mu}\de_0(m)d(x,y)\: \]
\item\label{odsSo1c}
$ k $ satisfies a linear decay condition:
	for each $ r\ge 0 \;$
	\[\ts
	\sup_{x\in Z}\sum_{m\in\Z}\vv{k(m,-(\cdot)\cdot_\tau x)\min\s{\v{m},r}}_\infty\le 2\:
	\]
\item\label{odsSo1d}
$ \ab{(\Z,\ell),\#} $ is strictly-polynomial-growing:
	for each $ r\ge 1 \;$
	\[ r\le \# B(r)\le 3r \:\]
\item\label{odsSo1e}
For all $ X,Y\in\Ss \;$
\nobreak\hfil\allowbreak\hfilneg
$
\smash{{r\ge 1\lRightarrow\dist_{\ssm H}\p{\spec\p{A_X},\spec\p{A_Y}}\le \rrm c\mspace{1mu minus 1mu}(\frac{r(r-1)}{2}+r+1)d_{\ssm H}(X,Y)+\frac*{2}{r}(\frac*{3^{4}}{2})^{\frac*{1}{2}}}}
$
and there exists $ \de\p{\rrm c}>0 $ (independent of $ X $ and $ Y $) such that
\[\ts
d_{\ssm H}(X,Y)\le{\min\s{\de(\rrm c),1}}
\lRightarrow
\dist_{\ssm H}\p{\spec\p{A_X},\spec\p{A_Y}}\le C d_{\ssm H}(X,Y)^{\frac*{1}{3}}\;
\]
where $ C\defeq \min\s{\rrm c\mspace{1mu minus 1mu}\p{\tfrac{\de^{\frac*{2}{3}}}{2}+\tfrac{\de^{\frac*{1}{3}}}{2}+1}+\tfrac{2}{\de^{\frac*{1}{3}}}\p{\tfrac{3^4}{2}}^{\frac*{1}{2}}:\de>0}\: $
Also, $ \de\p{\rrm c} $ minimizes $ C \:$
\end{enumerate}
\end{thm}
\noindent
Theorem~\ref{thm:odsSo1}(\ref{odsSo1a}, \ref{odsSo1b}, \ref{odsSo1c}, \ref{odsSo1d}) are easy to see and come in handy since Theorem~\ref{thm:odsSo1}(\ref{odsSo1e}) follows from a combination of (\ref{odsSo1a}--\ref{odsSo1d}) and Lemma~\ref{lem:po(3)} (where amenability twinkles) and Lemma~\ref{lem:po(5)}.
We remark that this operator family can be found within \cite{MP2001_BGSSkew,MP2002_BourgainSquareAMO1of2,MP2011_KrugerSkew} and this regularity result is new.
We forgo generalizing Theorem~\ref{H} (where the group action is instead assumed to be power Lipschitz).
\subsection{One-dimensional Discrete Limit-periodic Schrödinger Operators}\label{subsec:odlSo}
We demonstrate that our unified theory applies for this titular setting provided a lenient condition replaces condition \textit{(\ref{Ha1})} within Theorem~\ref{H}; see Theorem~\ref{thm:odlSo1}(\ref{odlSo1a}).

Let $ \ab{\Z,\#} $ be the one-dimensional lattice endowed with the counting measure.
Let $ \ell $ be the metric on $ \Z $ defined by
\[\ts
\ell(n,m)=\v{-n+m}
\]
for every $ n,m\in\Z \:$
Define the space $ \bar{Z}\defeq [0,1]^{\Z} \:$
Let $ d $ be the metric on $ \bar{Z} $ defined by
\[\ts
d(x,y)=\sum_{m\ge 0}\frac{1}{2^m}\max\s{\v{x_m-y_m},\v{x_{-m}-y_{-m}}}
\]
for every $ x,y\in\bar{Z} \:$
Let $ \ab{\bar{Z},(\Z,\tau)} $ be the dynamical system defined by
	\[\ts
	-n\cdot_{\tau} x=
	(m\mapsto x(n+m))\eqdef T^{-n}x
	\]
	for every $ n\in\Z\, x\in \bar{Z} \:$
Define the phase (sub)space $ Z\defeq \bar{Z}\cap\s{\textup{limit-periodic}} \;$ i.e., for each $ x\in\bar{Z} \;$ $ x\in Z $ if and only if there exists a sequence $ \s{y_j}\subseteq \bar{Z} $ such that $ \#\s{-n\cdot_\tau y_j:n\in\Z}<+\infty $ for every $ j $ and $ \vv{x-y_j}_\infty\to 0 $ as $ j\to +\infty \:$
Define, for $ v\in C(Z) \;$ the kernel $ k_v=k:\Z\times Z\to\C $ by
	\[\ts
	k(m,x)=\begin{cases}
	1&\iph m\in\s{-1,1}\\
	v(x)&\iph m=0\\
	0&\elss
	\end{cases}
	\]
	for every $ m\in\Z\, x\in Z \:$
Define the operator family $ \s{A_{v,x}}_{v\in C(Z),x\in Z}\subseteq \Ls\p{\ell^2(\Z)} $ by
	\[\ts
	(A_{x}\psi)(n)=\sum_{m\in\Z}k(-n+m,-n\cdot_\tau x)\psi(m)=\psi(n+1)+\psi(n-1)+v(T^{-n}x)\psi(n) \]
	for every $ v\in C(Z)\, x\in Z\, \psi\in \ell^2(\Z)\, n\in\Z \:$
\begin{thm}[log-Hölder regularity]\label{thm:odlSo1}
Assume the potential $ v:Z\to\R $ is Lipschitz continuous, i.e., there exists $ \rrm c_v=\rrm c>0 $ such that for all $ x,y\in Z \;$ $ \v{v(x)-v(y)}\le \rrm c\mspace{1mu minus 1mu}d(x,y) \:$
Then
\begin{enumerate}[(a)]
\item\label{odlSo1a}
$ \tau $ satisfies an exponential Lipschitz continuity condition:
	for each $ n\in \Z \;$
	\[\ts d(-n\cdot_\tau x,-n\cdot_\tau y)\le 2^{\v{n}}d(x,y)\: \]
\item\label{odlSo1b}
$ k $ satisfies a Lipschitz continuity condition:
	for each $ m\in \Z \;$
	\[\ts \v{k(m,x)-k(m,y)}\le \rrm c\mspace{1mu minus 1mu}\de_0(m)d(x,y)\: \]
\item\label{odlSo1c}
$ k $ satisfies a linear decay condition:
	for each $ r\ge 0 \;$
	\[\ts
	\sup_{x\in Z}\sum_{m\in\Z}\vv{k(m,-(\cdot)\cdot_\tau x)\min\s{\v{m},r}}_\infty\le 2\:
	\]
\item\label{odlSo1d}
$ \ab{(\Z,\ell),\#} $ is strictly-polynomial-growing:
	for each $ r\ge 1 \;$
	\[ r\le \# B(r)\le 3r \:\]
\item\label{odlSo1e}
For all $ X,Y\in\Ss \;$
\nobreak\hfil\allowbreak\hfilneg
$
\smash{{r\ge 1\lRightarrow
\dist_{\ssm H}\p{\spec\p{A_X},\spec\p{A_Y}}\le \rrm c\mspace{1mu minus 1mu}2^rd_{\ssm H}(X,Y)+\frac*{2}{r}(\frac*{3^{4}}{2})^{\frac*{1}{2}}
}}
$
and for each $ \rho\in (0,1) \;$ there exists $ \de\p{\rho}>0 $ (independent of $ X $ and $ Y $) such that
\[\ts\hspace*{\leftmargin minus \leftmargin}
d_{\ssm H}(X,Y)\le{\min\s{\de(\rho),1}}
\lRightarrow
\dist_{\ssm H}\p{\spec\p{A_X},\spec\p{A_Y}}\le C (\log_2(({\ss d_{\ssm H}(X,Y)^{\rho}})^{-1}))^{-1}\;
\]
where $ C\defeq \rrm c+\tfrac{2}{1}\p{\tfrac{3^4}{2}}^{\frac*{1}{2}}\: $
Also, $ \de\p{\rho}=\de $ satisfies an inequality: $ \de^{1-\rho}\le(\log_2(({\ss \de^{\rho}})^{-1}))^{-1} \:$
\end{enumerate}
\end{thm}
\noindent
Theorem~\ref{thm:odlSo1}(\ref{odlSo1a}, \ref{odlSo1b}, \ref{odlSo1c}, \ref{odlSo1d}) are easy to see and come in handy since Theorem~\ref{thm:odlSo1}(\ref{odlSo1e}) follows from a combination of (\ref{odlSo1a}--\ref{odlSo1d}) and Lemma~\ref{lem:po(3)} (where amenability twinkles) and Lemma~\ref{lem:po(5)}.
We remark that this operator family can be found within \cite{MP1981_ASLim,MP2009_ALim,MP2010_GLim} and this regularity result is new.
We forgo generalizing Theorem~\ref{H} (where the group action is instead assumed to be exponential Lipschitz).
\clearpage
\section{Preliminaries}\label{sec:p}
\subsection{Locally-compact Hausdorff Groups}\label{subsec:l-cHg}
We provide a terse exposition of
	topological groups.
The focus is on
	the existence of a translation-invariant sigma-finite measure and
	the existence of a translation-invariant proper metric.
The topics on unimodularity and amenability are principal axes.

Say $G$ is a~\emph{topological group} if
	$G$ is a group with a topology $\Tc\;$ where
	the group product and the group inverse are continuous.
Let $G$ be a topological group.
Assume $G$ is both Hausdorff and locally-compact.
Say $\lambda$ is a~\emph{Radon measure} on $G$ if
	$\lambda$ is a nonzero Borel measure on $G\;$ where
	$\lambda$ is outer regular for Borel sets and inner regular for open sets and finite for compact sets.
Say $\lambda$ is a~\emph{left Haar measure} on $G$ if
	$\lambda$ is a left-invariant Radon measure on $G\:$
\begin{thm}\label{thm:l-cHg1}
Let $G$ be a locally-compact Hausdorff group.
\begin{enumerate}[(a)]
\item\label{l-cHg1a}
There exists a left Haar measure $ \lambda $ on $ G \:$
Also, $ \lambda $ is unique modulo scaling.
\item\label{l-cHg1b}
If $ G $ is sigma-compact, then $ \lambda $ is sigma-finite and $ \lambda $ is inner regular for Borel sets.
\item\label{l-cHg1c}
If $ G $ is first-countable, then there exists a left-invariant $\Tc$-generating metric on $ G \:$
\item\label{l-cHg1d}
If $ G $ is second-countable, then there exists a left-invariant proper $\Tc$-generating metric on $ G \:$
\item\label{l-cHg1e}
$ G $ is sigma-compact and $ G $ is first-countable if and only if $ G $ is second-countable.
\end{enumerate}
\end{thm}
\noindent
The proof of Theorem~\ref{thm:l-cHg1}(\ref{l-cHg1a}, \ref{l-cHg1b}, \ref{l-cHg1c}, \ref{l-cHg1d})
	can be found within \cite{GR2016_Folland,RA2014_Yeh,GR1955_MontgomeryZippin,GR1974_Struble}.
Theorem~\ref{thm:l-cHg1}(\ref{l-cHg1e}) follows from
	the observation that every sigma-compact space is a Lindelöf space and
	the observation that, for metric spaces, Lindelöfness and separability and second-countability are equivalent;
	see \cite{TOP1970_Willard}.

\subsubsection{Unimodularity}\label{subsubsec:u}
Let $ \ab{G,\lambda} $ be a second-countable locally-compact Hausdorff group endowed with a left Haar measure.
Let $ \Bc $ be the collection of all Borel subsets of $ G \:$
Define, for $ g\in G$ and $ B\in\Bc \;$
	$ \lambda^g(B)\defeq \lambda(Bg)$ and $ \lambda^*(B)\defeq \lambda(\inv{B}) \:$
Observe $ \lambda^g $ is a left Haar measure on $ G $ and there exists $ c(g)>0 $ such that $ \lambda^g=c(g)\lambda \:$
It can be shown that $ c:G\to (0,+\infty) \;$ called the modular function, is a continuous homomorphism.
Say $ \ab{G,\lambda} $ is \emph{unimodular} if
	$ c $ is identically $ 1 $ or, equivalently, $ \lambda^g=\lambda $ for every $ g\in G \:$
It can be shown that if $ \ab{G,\lambda} $ is unimodular, then $ \lambda^*=\lambda \:$

Observe if $ G $ is discrete or $ G $ is abelian or $ G $ is compact, then $ \ab{G,\lambda} $ is unimodular.
To see how unimodularity follows from $G/[G,G]$ being compact,
	the image of $ [G,G]\defeq \overline{\s{[g,h]:g,h\in G}} $ under $ c $ is trivial since $ c $ is a homomorphism and
	$ c(G)=c(G/[G,G]) $ is compact since $ c $ is continuous and
	the only compact subgroup of $ (0,+\infty) $ is $ \s{1} \:$
The details of this subsubsection and more examples can be found within \cite{GR2016_Folland} such as connected Lie groups that are nilpotent or semisimple.

\subsubsection{Amenability}\label{subsubsec:a}
Let $ \ab{G,\lambda} $ be a second-countable locally-compact Hausdorff group endowed with~a left Haar measure.
Say $ \mfr $ is a \emph{mean} on $ G $ if
	$ \mfr $ is a nonzero bounded linear functional on $ L^\infty(G) \;$ where
	$ \mfr(\vp)\ge 0 $ for every $ \vp\in L_+^\infty(G) $ and $ \mfr(\1_G)=1 \:$
Say $ \ab{G,\lambda} $ is \emph{amenable} if
	a left-invariant mean on $ G $ exists or, equivalently, for each $ \ep>0 $ and for each nonempty compact subset $ F $ of $ G \;$ there exists a compactly-supported continuous nonzero $ \chi\in L_+^2(G) $ such that $ \sup_{h\in F}\vv{\cramped{\chi^h}-\chi}< \ep\vv{\chi}\: $
Here $ {\cramped{\chi^h} : G \to [0,+\infty) : g \mapsto \chi(gh)} \:$

It can be shown that if $ G $ is abelian or $ G $ is compact, then $ \ab{G,\lambda} $ is amenable.
Also, if $ G $ is solvable, then $ \ab{G,\lambda} $ is amenable since $ G $ is an extension of an amenable group by an amenable closed normal subgroup.
We remark that if $ G $ has a closed subgroup isomorphic to the free group with two generators $ \rrm F_2 \;$ then $ \ab{G,\lambda} $ is not amenable.
There is an extent to which nonamenability is related to the presence of $ \rrm F_2 $; see \texttt{chapter 3} within \cite{GR1988_Paterson}:
	if $ G $ is almost-connected, then $ \ab{G,\lambda} $ is amenable if and only if $ G $ does not have a closed subgroup isomorphic to $ \rrm F_2 $;
	if $ G $ is a connected Lie group, then $ G $ is solvable (thus amenable) if and only if $ G $ does not have a closed subgroup isomorphic to $ \rrm F_2 \:$
The details of this subsubsection and more examples can be found within \cite{GR1984_Pier,GR1988_Paterson} such as groups that are uniformly-distributed or polynomial-growing.
\subsection{Dynamically-defined Operator Families}\label{subsec:d-dof}
We provide a terse exposition of
	operator families induced by elements of a $ C^* $-algebra itself induced by a topological dynamical system.
The focus is on their covariant representation and properties, e.g., uniform boundedness and strong continuity.

Let $ \ab{G,\lambda} $ be a second-countable locally-compact Hausdorff group endowed with a left Haar measure.
Define, for $ \vp\in L^1(G) \;$
	$ \int \vp(h)dh\defeq \int \vp d\lambda \:$
Let $ \ab{Z,(G,\alpha)} $ be a compact Hausdorff space endowed with a left continuous action.
Specifically, $ \alpha $ is a left continuous action if and only if $ \alpha:G\times Z\to Z $ is continuous and $ \alpha(e,x)=x $ for every $ x\in Z $ and $ \alpha(g,\alpha(h,x))=\alpha(gh,x) $ for every $ g,h\in G\, x\in Z \:$
Define, for $ g\in G$ and $ x\in Z \;$
	$ \alpha_gx\defeq \alpha(g,x) \:$
Define, for $ g\in G$ and $ f\in C(Z)$ and $ x\in Z \;$
	$ (\alpha_g f)(x)\defeq f(\alpha_{\inv{g}}x) \:$
Here $ C(Z) $ $(\text{or }C_b(Z))$ is the unital $ C^* $-algebra of all (bounded) continuous complex-valued functions on $ Z \:$
Define, for $ x\in Z \;$
	$ [\thickspace]_x:C_c(G,C(Z))\to\Ls\p{L^2(G)} $ by
\[\ts
[f_{(\cdot)}]_x=\int d\lambda(h)[(\alpha_{(\cdot)}f_h)(x)]R_h =(f_{\inv{g}h}(\alpha_{\inv{g}}x))
\]
for every $ f_{(\cdot)}\in C_c(G,C(Z)) \:$
Here $ C_c(G,C(Z)) $ is the linear space of all compactly-supported continuous $ C(Z) $-valued functions on $ G \:$
Also, $ [(\alpha_{(\cdot)}f_h)(x)]:L^2(G)\to L^2(G) \;$ a multiplication operator, is defined by $ ([(\alpha_{(\cdot)}f_h)(x)]\psi)(g)=f_h(\alpha_{\inv{g}}x)\psi(g) $ for every $ \psi\in L^2(G)\, g\in G \:$
Also, $ R_h:L^2(G)\to L^2(G) \;$ called the right regular representation, is defined by $ (R_h\psi)(g)=\psi(gh) $ for every $ \psi\in L^2(G)\, g\in G \:$
Also, $ (f_{\inv{g}h}(\alpha_{\inv{g}}x)) $ is a matrix whose $ (g,h) $-th component is $ f_{\inv{g}h}(\alpha_{\inv{g}}x) \:$
More conventionally,
\[\ts
([f_{(\cdot)}]_x\psi)(g)=\int f_{h}(\alpha_{\inv{g}}x)\psi(gh) dh=\int f_{\inv{g}h}(\alpha_{\inv{g}}x)\psi(h)dh
\]
for every $ x\in Z\, f_{(\cdot)}\in C_c(G,C(Z))\, \psi\in L^2(G)\, g\in G \:$
Observe $ [\thickspace]_x $ is linear and $ [f_{(\cdot)}]_x $ is linear.
Also, $ U_{\inv{j}}[(\alpha_{(\cdot)}f_h)(x)]U_{j}=[(\alpha_{(\cdot)}f_h)(\alpha_{\inv{j}}x)] $ and $ [U_j,R_h]=0 \:$
Here $ U_j:L^2(G)\to L^2(G) \;$ called the left regular representation, is defined by
$
(U_j\psi)(g)=\psi(\inv{j}g)
$
for every $ \psi\in L^2(G)\, g\in G \:$

We remark that $ C_c(G,C(Z)) $ is a ${}^*$-algebra and $ \vv{f_{(\cdot)}*\tilde{f}_{(\cdot)}}\le \vv{f_{(\cdot)}}\vv{\tilde{f}_{(\cdot)}} $ and $ \vv{f_{(\cdot)}\!{}^**f_{(\cdot)}}=\vv{f_{(\cdot)}}^2 \;$ where
$(f_{(\cdot)}*\tilde{f}_{(\cdot)})_h(x)\defeq \int f_j(x)\tilde{f}_{\inv{j}h}(\alpha_{\inv{j}}x)dj\, (f_{(\cdot)}\!{}^*)_h(x)\defeq \overline{f_{\inv{h}}(\alpha_{\inv{h}}x)}\, \vv{f_{(\cdot)}}\defeq\sup_{x\in Z}\vv{[f_{(\cdot)}]_x}\:$
Also, $ C_c(G,C(Z)) $ embeds into a $ C^* $-algebra, denoted $ C(Z)\rtimes_{\alpha} G \;$ where
\[\ts \spec(f_{(\cdot)})=\overline{\bigcup_{x\in Z}\spec([f_{(\cdot)}]_x)} \]
for normal elements; see \cite{MP1985_BellissardLimaTestard,MP1992_BellissardGLT,LSA2018_Pedersen,MP2018_BeckusBellissardDeNittisGCT1of2}.
The $ C^* $-algebra $C(Z)\rtimes_{\alpha} G\;$ called a reduced crossed product, and its $ K $-theory were utilized to prove the Gap Labeling Theorems; see \cite{MP1985_BellissardLimaTestard,MP1986_Bellissard,MP1992_BellissardBovierGhezGLT,MP1992_BellissardGLT,GR2001_Lindenstrauss}:
if $ \ab{G,\lambda} $ is both unimodular and amenable and $Z$ is second-countable, then for each ergodic left-invariant Borel probability measure $\mu$ on $ Z \;$ there exists a countable subgroup of $ \R $ containing all gap labels.
The 1985 result by J.~Bellissard; R.~Lima; D.~Testard \cite{MP1985_BellissardLimaTestard} and the 1986 result by J.~Bellissard \cite{MP1986_Bellissard} considered the case where $ G=\Z^d $ or $ G=\R^d $ or $ G $ is both discrete and amenable and utilized that $G$ sometimes has a regular Følner sequence so that the pointwise ergodic theorem holds for the `time-average' approximations of the `space-average' integrated density of states---whose image over the gaps in the spectra of self-adjoint elements of $ C(Z)\rtimes_{\alpha} G $ are precisely the gap labels.
The existence of regular Følner sequences is not known for amenable groups with exponential growth \cite{GR2006_Nevo}, but the 2001 result by E.~Lindenstrauss \cite{GR2001_Lindenstrauss} establishes that amenable groups always have tempered Følner sequences.
We remark that tempered Følner sequences are suitable.

This subsection is the only place where $ C^* $-algebras are mentioned, and it is done so to expose the spectrum map
\[\ts X\mapsto \overline{\bigcup_{x\in X}\spec([f_{(\cdot)}]_x)}\;\] where $ X $ is a subsystem of $ \ab{Z,(G,\alpha)} \:$
Specifically, $ X $ is a subsystem of $ \ab{Z,(G,\alpha)} $ if and only if $ X $~is a left-invariant nonempty closed subset of $ Z \:$

\begin{lem}\label{lem:d-dof1}
Let $ \ab{G,\lambda} $ be a second-countable locally-compact Hausdorff group endowed with a left Haar measure.
Let $ \ab{Z,(G,\alpha)} $ be a topological space endowed with a left continuous action.
Let $ f_{(\cdot)}\in C_c(G,C_b(Z)) \:$
Define, for $ x\in Z \;$
	${A}_x:L^2(G)\to L^2(G)$ by
\[\ts
({A}_x\psi)(g)=\int f_{h}(\alpha_{\inv{g}}x)\psi(gh) dh=\int f_{\inv{g}h}(\alpha_{\inv{g}}x)\psi(h)dh
\]
for every $ \psi\in L^2(G)\, g\in G \:$
Define, for $ (h,x)\in G\times Z\, g_{\ttm 1},g_{\ttm 2}\in G \;$
	$ k(\inv{g_{\ttm 1}}h,\inv{g_{\ttm 2}}x)\defeq f_{\inv{g_{\ttm 1}}h}(\alpha_{\inv{g_{\ttm 2}}}x) \:$
Let $ \Ss $ be the collection of all subsystems of $ \ab{Z,(G,\alpha)} \:$
Define, for $ X\in\Ss \;$
\[\ts \spec\p{A_X}\defeq\overline{\bigcup_{x\in X}\spec\p{A_x}} \:\]
\begin{enumerate}[(a)]
\item\label{d-dof1a}
$ f_h(x)=k(h,x)=k:G\times Z\to \C $ is measurable and $ \sup_{x\in Z}\int\vv{k(h,\inv{(\cdot)}x)}_{\infty}dh<+\infty \:$
\item\label{d-dof1b}
If $ \ab{G,\lambda} $ is unimodular, then for all $ x,y\in Z $ and for each $ \psi\in L^2(G) $ and for each $ j\in G \;$
\begin{align*}
\ts\vv{{A}_x\psi}\le& (\int
\vv{k(h,\inv{(\cdot)}x)\1_{\supp(\psi)\inv{h}}(\cdot)}_{\infty}
dh)\vv{\psi}\;\\
\ts\vv{{A}_x\!{}^{*}\psi}\le& (\int
\vv{k(h,\inv{(\cdot)}x)\1_{\supp(\psi)}(\cdot)}_{\infty}
dh)\vv{\psi}\;\\
\ts\vv{({A}_y-{A}_x)\psi}\le& (\int
\vv{\p{k(h,\inv{(\cdot)}y)-k(h,\inv{(\cdot)}x)}\1_{\supp(\psi)\inv{h}}(\cdot)}_{\infty}
dh)\vv{\psi}\;\\
\ts\vv{({A}_y-{A}_x)^{*}\psi}\le& (\int
\vv{\p{k(h,\inv{(\cdot)}y)-k(h,\inv{(\cdot)}x)}\1_{\supp(\psi)}(\cdot)}_{\infty}
dh)\vv{\psi}\;\\
\ts U_{\inv{j}}{A}_xU_{j}=&{A}_{\inv{j}x}\:
\end{align*}
\item\label{d-dof1c}
If $ \ab{G,\lambda} $ is unimodular and $ A_x $ is normal for every $ x\in Z \;$ then for each $ x\in Z \;$
\[\ts \lim_{y\xrightarrow[]{\smash{Z}}x}\sup_{E\in\spec({A}_x)}\dist(E,\spec({A}_y))=0 \]
and $\spec\p{A_x}=\spec\p{A_{X(x)}}\;$ where $ X\p{x}\defeq \overline*{\s{\inv{g}x:g\in G}} \:$
\end{enumerate}
\end{lem}
\noindent
Lemma~\ref{lem:d-dof1}(\ref{d-dof1a}) follows from
	$ f_{(\cdot)}\in C_c\p{G,C_b(Z)} $ and $ f_h(x)=k(h,x) \:$
The proof of---a suitable generalization of---Lemma~\ref{lem:d-dof1}(\ref{d-dof1b}, \ref{d-dof1c}) can be found within \textsc{section~\ref{sec:A1}}.
\subsection{Hausdorff Distance}\label{subsec:hd}
Let $ Z $ be a metrizable space.
Let $ d $ be a $ \Tc $-generating metric on $ Z \:$
Specifically, for each metric $ d $ on $ Z \;$
	$ d $ is $ \Tc $-generating if and only if the balls $ \s{y\in Z:d(y,x)<r} $ are open and generate the topology on $ Z \:$
	Also, the balls are precompact if and only if $ d $ is proper.
Let $ S_{\ttm 1} $ and $ S_{\ttm 2} $ be subsets of $ Z \:$
The \textit{Hausdorff distance} between $ S_{\ttm 1} $ and $ S_{\ttm 2} $ with respect to $ d $ is
\[\ts \max\s{\sup_{s\in S_{\ttm 1}}d(s,S_{\ttm 2}),\sup_{s\in S_{\ttm 2}}d(s,S_{\ttm 1})}\eqdef d_{\mathsf{H}}(S_{\ttm 1},S_{\ttm 2}) \:\]
\begin{thm}\label{thm:hd1}
Let $ G $ be a second-countable locally-compact Hausdorff group.
Let $ \ab{Z,(G,\alpha)} $ be a compact metrizable (or second-countable compact Hausdorff) space endowed with a left continuous action.
Let $ \Ss $ be the collection of all subsystems of $ \ab{Z,(G,\alpha)} \:$
Let $ d $ be a $ \Tc $-generating metric on $ Z \:$
Then $ (\Ss,d_{\ssm H}) $ is a compact metric space.
\end{thm}
\noindent
The proof of Theorem~\ref{thm:hd1} can be found within \cite{MP2018_BeckusBellissardDeNittisGCT1of2}.
\subsection{Strict Polynomial Growth}\label{subsec:spg}
Let $ \ab{G,\lambda} $ be a second-countable locally-compact Hausdorff group endowed with a left Haar measure.
Let $ \ell $ be a left-invariant proper $ \Tc $-generating metric~on~$ G \:$
Define, for $ g\in G \;$ $ \v{g}\defeq \ell(g,e) \:$
Define, for $ r\ge 0 \;$ $ B(r)\defeq \s{g\in G:\v{g}<r} \:$
Say $ \ab{(G,\ell),\lambda} $ is \textit{strictly-polynomial-growing} if there exist $ b\ge 1 $ and $ c_1\ge c_0>0 $ such that for each $ r\ge 1 \;$
\[\ts c_0r^b\le \lambda(B(r))\le c_1r^b \:\]
\begin{thm}\label{thm:spg1}
Let $ \ab{G,\lambda} $ be a second-countable locally-compact Hausdorff group endowed with a left Haar measure.
Let $ \ell $ be a left-invariant proper $ \Tc $-generating metric on $ G \:$
Assume $ \ab{(G,\ell),\lambda} $ is strictly-polynomial-growing.
Then $ \ab{G,\lambda} $ is both unimodular and amenable.
\end{thm}
\noindent
The proof of Theorem~\ref{thm:spg1} can be found within \cite{GR2006_Nevo}.
\clearpage
\section{Proof of Main Results}\label{sec:po}
{
\noindent
Let $ \ab{G,\lambda} $ be
	a second-countable locally-compact Hausdorff group endowed with
	a left Haar measure.
Let $ \ab{Z,(G,\alpha)} $ be
	a metrizable space endowed with
	a left continuous action.
Let $ k:G\times Z\to\C \;$
	where $ k $ is measurable and
	$\sup_{x\in Z}\int \vv{k(h,\inv{(\cdot)}x)}_{\infty} dh<+\infty\:$
Define, for $ x\in Z \;$
	${A}_x:L^2(G)\to L^2(G)$ by
\[\ts
(A_x\psi)(g)=\int k(h,\inv{g}x)\psi(gh) dh=\int k(\inv{g}h,\inv{g}x)\psi(h)dh
\]
for every $ \psi\in L^2(G)\, g\in G \:$
Let $ \Ss $ be the collection of all subsystems of $ \ab{Z,(G,\alpha)} \:$
Define, for $ X\in\Ss \;$
	$ \spec(A_X)\defeq \overline{\bigcup_{x\in X}\spec(A_x)} \:$
Let $ \ell $ be a left-invariant proper $ \Tc $-generating metric on $ G \:$
Let $ d $ be a $ \Tc $-generating metric on $ Z \:$
}
\begin{thm}\label{thm:po(1)}
Assume $ \ab{G,\lambda} $ is unimodular, $ A_x $ is normal for every $ x\in Z \;$ and the following.
\begin{enumerate}[(i)]
\item
$ \alpha $ satisfies a uniform continuity condition:
	for each $ \ep>0 $ and for each nonempty compact subset $ K $ of $ G \;$ there exists $ \de>0 $ such that
	$ d(y,x)<\de \lRightarrow \sup_{g\in K}d(\inv{g}y,\inv{g}x)<\ep \:$
\item
$ k $ satisfies a uniform continuity condition:
	for each $ \ep>0 $ and for each nonempty compact subset $ F $ of $ G \;$ there exists $ \de>0 $ such that
	$ d(y,x)<\de \lRightarrow \int\v{k(h,y)-k(h,x)}\1_F\p{h}dh<\ep\:$
\item
$ k $ satisfies a uniform decay condition:
	for each $ \ep>0 \;$ there exists a nonempty compact subset $ F $ of $ G $ such that
	$ \sup_{x\in Z}\int\vv{k(h,\inv{(\cdot)}x)\1_{G\setminus F}(h)}_{\infty}dh<\ep \:$
\end{enumerate}
Fix $ X\in\Ss \:$
Then for each $ \ep>0 \;$ there exists $ \de>0 $ such that for each $ Y\in\Ss \;$ if $ d_{\ssm H}(X,Y)<\de \;$ then $ \sup_{E\in\spec\p{A_X}}\dist(E,\spec\p{A_Y})<\ep \;$ i.e.,
\[\ts \lim_{Y\xrightarrow[]{\smash{\Ss}}X}\sup_{E\in\spec(A_X)}\dist(E,\spec(A_Y))=0\: \]
\end{thm}
\begin{proof}
It suffices to show that
	for each $ E\in\spec\p{A_X} $ and for each $ \ep>0 \;$ there exists $ \de>0 $ such that for each $ Y\in \Ss \;$ if $ d_{\ssm H}(X,Y)<\de \;$ then $ \dist(E,\spec\p{A_Y})<\ep $ since $ \spec\p{A_X} $ is compact.
Fix $ E\in\spec\p{A_X} $ and $ \ep>0 \:$
Observe there exists a nonempty compact subset $ F $ of $ G $ such that
\[\ts
\sup_{x\in Z}\int\vv{k(h,\inv{(\cdot)}x)\1_{G\setminus F}(h)}_{\infty}dh<\frac{1}{5}\ep
\]
since $ k $ satisfies a uniform decay condition.
Observe there exist $ x\in X $ and $ E'\in\spec\p{A_x} $ such that $ \v{E-E'}<\frac{1}{5}\ep \:$
Without loss of generality, assume $ E'\notin\spec\p{A_Y} \:$
Observe $ A_x $ is normal and, by Weyl's criterion, there exists a nonzero vector $ \vp $ such that $ \vv{(A_x-E')\vp}<\frac{1}{5}\ep\vv{\vp} \:$
Without loss of generality, assume $ \vp $ is both continuous and compactly-supported.
Define $ K\defeq\supp\p{\vp} \:$
Observe there exists $ \de>0 $ such that
\[\ts
d(y,x)<2\de\lRightarrow
\int\vv{\p{k(h,\inv{(\cdot)}y)-k(h,\inv{(\cdot)}x)}\1_K(\cdot)}_\infty\1_F(h)dh<\frac{1}{5}\ep
\]
since the group action and $ k $ satisfy a uniform continuity condition.
Fix $ Y\in\Ss \:$
Assume $ d_{\ssm H}(X,Y)<\de \:$
Observe there exists $ y\in Y $ such that $d(x,y)< d(x,Y)+\de<2\de\:$
Also, $ A_y $ is normal and
\[\ts \vv{(A_y-E')^{-1}}^{-1}\vv{\vp}\le\vv{(A_y-E')\vp}=\vv{(A_y-E')^*\vp}<\p{\vv{(A_y-A_x)^*\frac{\vp}{\vv{\vp}}}+\frac{1}{5}\ep}\vv{\vp}\]
and
\begin{align*}\swapabovedisplayskip
&\dist\p{E,\spec\p{A_Y}}\\
&< \dist\p{E',\spec\p{A_y}}+\frac{1}{5}\ep\\
&=\vv{(A_y-E')^{-1}}^{-1}+\frac{1}{5}\ep\\
&<\vv{(A_y-A_x)^*\frac{\vp}{\vv{\vp}}}+\frac{2}{5}\ep\\
&\le\int\vv{(k(h,\inv{(\cdot)}y)-k(h,\inv{(\cdot)}x))\1_{K}(\cdot)}_{\infty}dh+\frac{2}{5}\ep
\tagps{\textrm{Prop~\ref{prop:A1(1)}}}\\
&<\int\vv{(k(h,\inv{(\cdot)}y)-k(h,\inv{(\cdot)}x))\1_{K}(\cdot)}_{\infty}\1_{F}(h)dh+\frac{4}{5}\ep\\
&< \ep\:\qedhere
\end{align*}
\end{proof}
\begin{rmk}\label{rmk:po(1.1)}
Theorem~\ref{thm:po(1)} will be improved provided the group is amenable.
\end{rmk}
\noindent
{
Let $ \chi $ be a compactly-supported continuous nonzero vector in $ L^2\p{G} \:$
The \emph{cutoff operator} with respect to $ \chi $ is $ [\chi]:L^2(G)\to L^2(G) $ defined by $ ([\chi]\psi)(g)=\chi(g)\psi(g) $ for every $ \psi\in L^2(G)\, g\in G \:$
Let $ A $ and $ B $ be operators on $ L^2(G) \:$
The \emph{commutator operator} with respect to $ A $ and $ B $ is $ [A,B]:L^2(G)\to L^2(G) $ defined by $ ([A,B]\psi)(g)=((AB-BA)\psi)(g) $ for every $ \psi\in L^2(G)\, g\in G \:$
}
\begin{lem}\label{lem:po(2)}
\newcommand{\chih}{\cramped{\chi^h}}%
Assume $ \ab{G,\lambda} $ is unimodular.
Let $ \chi $ be a compactly-supported continuous nonzero vector.
Define, for $ j,h\in G \;$
\begin{align*}\swapabovedisplayskip
\chi_j:&G\to\C:g\mapsto\chi\p{\inv{j}g}\;\\
\chih:&G\to\C:g\mapsto\chi\p{gh}\:
\end{align*}
Fix $ x\in Z$ and $ \psi\in L^2(G) \:$
Then
\begin{align*}
\int\vv{[A_x,[\chi_j]]\psi}^2dj
	\le&\p{\int \vv{k(h,\inv{(\cdot)}x)\1_{\supp\p{\psi}\inv{h}}(\cdot)}_{\infty}\vv{\chih-\chi}dh}^2\vv{\psi}^2\;\\
\int\vv{[A_x\!{}^*,[\chi_j]]\psi}^2dj
	\le&\p{\int \vv{k(\inv{h},\inv{(\cdot)}x)\1_{\supp\p{\psi}}(\cdot)}_{\infty}\vv{\chih-\chi}dh}^2\vv{\psi}^2\;\\
\int\vv{[\chi_j]\psi}^2dj=&\vv{\chi}^2\vv{\psi}^2\:
\end{align*}
\end{lem}
\begin{proof}
\newcommand{\chih}{\cramped{\chi^h}}%
Observe
\begin{align*}
&\int\vv{[A_x,[\chi_j]]\psi}^2dj\\
&=\iint\v{\int k(h,\inv{g}x)(\chi_j(gh)-\chi_j(g))\psi(gh)dh}^2dgdj\\
&\le\iint\p{\int \v{k(h,\inv{g}x)}\v{\chi_j(gh)-\chi_j(g)}\v{\psi(gh)}dh}^2dgdj\\
&\le\iint\!\begin{array}[t]{l}
\p{\int \v{k(h,\inv{g}x)}\1_{\supp\p{\psi}}(gh)\vv{\chih-\chi}dh}\\
\p{\int \v{k(h,\inv{g}x)}\v{\chi_j(gh)-\chi_j(g)}^2\v{\psi(gh)}^2\vv{\chih-\chi}^{-1}dh}dgdj
\end{array}
\tagps{\textrm{Hölder's}\\\textrm{inequality}}\\
&\le\iint\!\begin{array}[t]{l}
\p{\int \vv{k(h,\inv{(\cdot)}x)\1_{\supp\p{\psi}\inv{h}}(\cdot)}_{\infty}\vv{\chih-\chi}dh}\\
\p{\int \vv{k(h,\inv{(\cdot)}x)\1_{\supp\p{\psi}\inv{h}}(\cdot)}_{\infty}\v{\chi_j(gh)-\chi_j(g)}^2\v{\psi(gh)}^2\vv{\chih-\chi}^{-1}dh}djdg
\end{array}
\\
&=\p{\int \vv{k(h,\inv{(\cdot)}x)\1_{\supp\p{\psi}\inv{h}}(\cdot)}_{\infty}\vv{\chih-\chi}dh}^2\vv{\psi}^2
\shortintertext{and}
&\int\vv{[A_x\!{}^*,[\chi_j]]\psi}^2dj\\
&=\iint\v{\int \overline{k(\inv{h},\inv{h}\inv{g}x)}(\chi_j(gh)-\chi_j(g))\psi(gh)dh}^2dgdj\\
&\le\iint\p{\int \v{\overline{k(\inv{h},\inv{h}\inv{g}x)}}\v{\chi_j(gh)-\chi_j(g)}\v{\psi(gh)}dh}^2dgdj\\
&\le\iint\!\begin{array}[t]{l}
\p{\int \v{k(\inv{h},\inv{(gh)}x)}\1_{\supp\p{\psi}}(gh)\vv{\chih-\chi}dh}\\
\p{\int \v{k(\inv{h},\inv{(gh)}x)}\v{\chi_j(gh)-\chi_j(g)}^2\v{\psi(gh)}^2\vv{\chih-\chi}^{-1}dh}dgdj
\end{array}
\\
&\le\iint\!\begin{array}[t]{l}
\p{\int \vv{k(\inv{h},\inv{(\cdot)}x)\1_{\supp\p{\psi}}(\cdot)}_{\infty}\vv{\chih-\chi}dh}\\
\p{\int \vv{k(\inv{h},\inv{(\cdot)}x)\1_{\supp\p{\psi}}(\cdot)}_{\infty}\v{\chi_j(gh)-\chi_j(g)}^2\v{\psi(gh)}^2\vv{\chih-\chi}^{-1}dh}djdg
\end{array}
\\
&=\p{\int \vv{k(\inv{h},\inv{(\cdot)}x)\1_{\supp\p{\psi}}(\cdot)}_{\infty}\vv{\chih-\chi}dh}^2\vv{\psi}^2\:
\end{align*}
Also,
\[\ts
\int\vv{[\chi_j]\psi}^2dj
=\iint\v{\chi_j(g)}^2\v{\psi(g)}^2dgdj
=\iint\v{\chi_j(g)}^2\v{\psi(g)}^2djdg
=\vv{\chi}^2\vv{\psi}^2\:\qedhere
\]
\end{proof}
\begin{lem}\label{lem:po(3)}
\newcommand{\chih}{\cramped{\chi^h}}%
Assume $ \ab{G,\lambda} $ is unimodular and $ A_x $ is normal for every $ x\in Z \:$
Let $ \chi $ be a compactly-supported continuous nonzero vector.
Fix $ X,Y\in\Ss \:$
Then for each $ \rho>0 \;$
\[\ts
\Delta_{\ssm H}
\le\sup_{\lsubstack{\check{x}\in X,\check{y}\in Y,\check{z}\in Z\\\mathrlap{d(\check{x},\check{y})< \de_{\ssm H}+\rho}}}
	\int\vv{(k(h,\inv{(\cdot)}\check{y})-k(h,\inv{(\cdot)}\check{x}))\1_{\supp\p{\chi}}(\cdot)}_{\infty}+ \vv{k(h,\inv{(\cdot)}\check{z})}_{\infty}\frac{\vv{\chih-\chi}}{\vv{\chi}}dh\;
\]
where $ \Delta_{\ssm H}\defeq \dist_{\ssm H}\p{\spec\p{A_X},\spec\p{A_Y}}$ and $ \de_{\ssm H}\defeq d_{\ssm H}(X,Y)$ and $ \chih:G\to\C:g\mapsto\chi\p{gh} \:$
\end{lem}
\begin{proof}
\newcommand{\chih}{\cramped{\chi^h}}%
Fix $ \rho >0 \:$
Let $ E\in\spec\p{A_X} \:$
Let $ \ep>0 \:$
Observe there exist $ x\in X$ and $ E'\in\spec\p{A_x} $ such that $ \v{E-E'}<\ep \:$
Without loss of generality, assume $ E'\notin\spec\p{A_Y} \:$
Observe $ A_x $ is normal and, by Weyl's criterion, there exists a nonzero vector $ \vp $ such that $\vv{(A_x-E')\vp}<\ep\vv{\vp}\:$
Without loss of generality, assume $ \vp $ is both continuous and compactly-supported.
Define, for $ j\in G \;$
\[\ts \chi_j:G\to\C:g\mapsto\chi\p{\inv{j}g}\:\]
Define $ r\defeq \frac{1}{\ep}\int \vv{k(h,\inv{(\cdot)}x)}_{\infty}\frac{\vv{\chih-\chi}}{\vv{\chi}}dh \:$
Observe
\begin{align*}
&\int\vv{(A_x-E')[\chi_j]\vp}^2dj\\
&=\int\vv{[(A_x-E'),[\chi_j]]\vp+[\chi_j](A_x-E')\vp}^2dj\\
&\le \int  (1+\frac{1}{r})\vv{[(A_x-E'),[\chi_j]]\vp}^2+(1+r)\vv{[\chi_j](A_x-E')\vp}^2dj
\tagps{\textrm{Peter--Paul's}\\\textrm{inequality}}\\
&\le (1+\frac{1}{r})\p{\int \vv{k(h,\inv{(\cdot)}x)}_{\infty}\vv{\chih-\chi}dh}^2\vv{\vp}^2+(1+r)\vv{\chi}^2\vv{(A_x-E')\vp}^2
\tagps{\textrm{Lem~\ref{lem:po(2)}}}\\
&< (1+\frac{1}{r})\p{\int \vv{k(h,\inv{(\cdot)}x)}_{\infty}\vv{\chih-\chi}dh}^2\vv{\vp}^2+(1+r)\ep^2\vv{\chi}^2\vv{\vp}^2\\
&=\int
(1+\frac{1}{r})\p{\int \vv{k(h,\inv{(\cdot)}x)}_{\infty}\vv{\chih-\chi}dh}^2\vv{\chi}^{-2}\vv{[\chi_j]\vp}^2+(1+r)\ep^2\vv{[\chi_j]\vp}^2
dj\\
&=\int
\Bigp{(1+\frac{1}{r})\p{\int \vv{k(h,\inv{(\cdot)}x)}_{\infty}\frac{\vv{\chih-\chi}}{\vv{\chi}}dh}^2+(1+r)\ep^2}\vv{[\chi_j]\vp}^2
dj\\
&=\int
\p{\int \vv{k(h,\inv{(\cdot)}x)}_{\infty}\frac{\vv{\chih-\chi}}{\vv{\chi}}dh+\ep}^2\vv{[\chi_j]\vp}^2
dj\:
\end{align*}
Observe there exists $ j\in G $ such that
\[\ts\vv{(A_x-E')[\chi_j]\vp}<\p{\int \vv{k(h,\inv{(\cdot)}x)}_{\infty}\frac*{\vv{\chih-\chi}}{\vv{\chi}}dh+\ep}\vv{[\chi_j]\vp}\:\]
Define~$ \vp_j\defeq [\chi_j]\vp \:$
Observe
\[\ts
\supp\p{\vp_j}\subseteq\supp\p{\chi_j}\:
\]
Observe there exists $ y\in Y $ such that $d(\inv{j}x,\inv{j}y)< d(\inv{j}x,Y)+\rho\le \de_{\ssm H}+\rho$ since $Y$ is left-invariant.
Also, $ A_y $ is normal and
\[\ts
\vv{(A_y-E')\vp_j}
=\vv{(A_y-E')^{*}\vp_j}
<\p{\vv{(A_y-A_x)^*\frac{\vp_j}{\vv{\vp_j}}}+\int \vv{k(h,\inv{(\cdot)}x)}_{\infty}\frac{\vv{\chih-\chi}}{\vv{\chi}}dh+\ep}\vv{\vp_j}
\]
and
\begin{align*}
&\dist\p{E,\spec\p{A_Y}}\\
&<\dist\p{E',\spec\p{A_y}}+\ep\\
&=\vv{(A_y-E')^{-1}}^{-1}+\ep\\
&\le \vv{(A_y-E')\frac{\vp_j}{\vv{\vp_j}}}+\ep\\
&<\vv{(A_y-A_x)^*\frac{\vp_j}{\vv{\vp_j}}}+\int \vv{k(h,\inv{(\cdot)}x)}_{\infty}\frac{\vv{\chih-\chi}}{\vv{\chi}}dh+2\ep\\
&\le\int\vv{(k(h,\inv{(\cdot)}y)-k(h,\inv{(\cdot)}x))\1_{\supp\p{\vp_j}}(\cdot)}_{\infty}+ \vv{k(h,\inv{(\cdot)}x)}_{\infty}\frac{\vv{\chih-\chi}}{\vv{\chi}}dh+2\ep\\
&\le\int\vv{(k(h,\inv{(\cdot)}y)-k(h,\inv{(\cdot)}x))\1_{\supp\p{\chi_j}}(\cdot)}_{\infty}+ \vv{k(h,\inv{(\cdot)}x)}_{\infty}\frac{\vv{\chih-\chi}}{\vv{\chi}}dh+2\ep\\
&=\int\vv{(k(h,\inv{(\cdot)}\inv{j}y)-k(h,\inv{(\cdot)}\inv{j}x))\1_{\supp\p{\chi}}(\cdot)}_{\infty}+ \vv{k(h,\inv{(\cdot)}x)}_{\infty}\frac{\vv{\chih-\chi}}{\vv{\chi}}dh+2\ep\\
&\le\sup_{\lsubstack{\check{x}\in X,\check{y}\in Y,\check{z}\in Z\\\mathrlap{d(\check{x},\check{y})< \de_{\ssm H}+\rho}}}
	\int\vv{(k(h,\inv{(\cdot)}\check{y})-k(h,\inv{(\cdot)}\check{x}))\1_{\supp\p{\chi}}(\cdot)}_{\infty}+ \vv{k(h,\inv{(\cdot)}\check{z})}_{\infty}\frac{\vv{\chih-\chi}}{\vv{\chi}}dh+2\ep\:
\end{align*}
As a result, $ \sup_{E\in\spec\p{A_X}}\dist\p{E,\spec\p{A_Y}} $ and $ \Delta_{\ssm H} $ are bounded from above by
\[\ts
\sup_{\lsubstack{\check{x}\in X,\check{y}\in Y,\check{z}\in Z\\\mathrlap{d(\check{x},\check{y})< \de_{\ssm H}+\rho}}}
	\int\vv{(k(h,\inv{(\cdot)}\check{y})-k(h,\inv{(\cdot)}\check{x}))\1_{\supp\p{\chi}}(\cdot)}_{\infty}+ \vv{k(h,\inv{(\cdot)}\check{z})}_{\infty}\frac{\vv{\chih-\chi}}{\vv{\chi}}dh
\]
as desired.
\end{proof}
\begin{thm}\label{thm:po(4)}
Assume $ \ab{G,\lambda} $ is unimodular, $ A_x $ is normal for every $ x\in Z \;$ and the following.
\begin{enumerate}[(i)]
\item
$ \alpha $ satisfies a uniform continuity condition:
	for each $ \ep>0 $ and for each nonempty compact subset $ K $ of $ G \;$ there exists $ \de>0 $ such that
	$ d(y,x)<\de \lRightarrow \sup_{g\in K}d(\inv{g}y,\inv{g}x)<\ep \:$
\item
$ k $ satisfies a uniform continuity condition:
	for each $ \ep>0 $ and for each nonempty compact subset $ F $ of $ G \;$ there exists $ \de>0 $ such that
	$ d(y,x)<\de \lRightarrow {\int\v{k(h,y)-k(h,x)}\1_F\p{h}dh}<\ep\:$
\item
$ k $ satisfies a uniform decay condition:
	for each $ \ep>0 \;$ there exists a nonempty compact subset $ F $ of $ G $ such that
	$ \sup_{x\in Z}\int\vv{k(h,\inv{(\cdot)}x)\1_{G\setminus F}(h)}_{\infty}dh<\ep \:$
\item
$ \ab{G,\lambda} $ is amenable.
\end{enumerate}
Then for each $ \ep>0 \;$ there exists $ \de>0 $ such that for all $ X,Y\in\Ss \;$ if $ d_{\ssm H}(X,Y)<\de \;$ then $\dist_{\ssm H}(\spec\p{A_X},\spec\p{A_Y})<\ep\;$ i.e.,
\[\ts
\lim_{\de\to 0^+}\sup_{X,Y\in\Ss,d_{\ssm{H}}(X,Y)< \de}\dist_{\ssm{H}}(\spec(A_X),\spec(A_Y))=0\:
\]
\end{thm}
\begin{proof}
\newcommand{\chih}{\cramped{\chi^h}}%
Fix $ \ep>0 \:$
Observe there exists a nonempty compact subset $ F $ of $ G $ such that
\[\ts
\sup_{x\in Z}\int\vv{k(h,\inv{(\cdot)}x)\1_{G\setminus F}(h)}_{\infty}dh<\frac{1}{6}\ep
\]
since $ k $ satisfies a uniform decay condition.
Define $\vv{k}\defeq \sup_{x\in Z}\int \vv{k(h,\inv{(\cdot)}x)}_\infty dh \:$
Observe there exists a compactly-supported continuous nonzero vector $ \chi\in L^2_+(G) $ such that
\[\ts \sup_{h\in F}\frac{\vv{\chih-\chi}}{\vv{\chi}}<\frac{1}{6(\vv{k}+1)}\ep
 \]
since $ \ab{G,\lambda} $ is amenable; see \textsc{section~\ref{subsubsec:a}}.
Here $ \chih:G\to \lbrp{0,+\infty}:g\mapsto \chi(gh) \:$
Define $ K\defeq\supp\p{\chi} \:$
Observe there exists $ \de>0 $ such that
\[\ts
d(y,x)<2\de\lRightarrow
\int\vv{\p{k(h,\inv{(\cdot)}y)-k(h,\inv{(\cdot)}x)}\1_K(\cdot)}_\infty\1_F(h)dh<\frac{1}{6}\ep
\]
since the group action and $ k $ satisfy a uniform continuity condition.
Fix $ X,Y\in\Ss \:$
Assume $ d_{\ssm H}(X,Y)<\de \:$
Observe
\begin{align*}
&\dist_{\ssm H}\p{\spec\p{A_X},\spec\p{A_Y}}\\
&\le\sup_{\lsubstack{\check{x}\in X,\check{y}\in Y,\check{z}\in Z\\\mathrlap{d(\check{x},\check{y})< d_{\ssm H}(X,Y)+\de}}}
	\int\vv{(k(h,\inv{(\cdot)}\check{y})-k(h,\inv{(\cdot)}\check{x}))\1_{K}(\cdot)}_{\infty}+ \vv{k(h,\inv{(\cdot)}\check{z})}_{\infty}\frac*{\vv{\chih-\chi}}{\vv{\chi}}dh
\tags{\p{\lsubstack{\textrm{Lem~\ref{lem:po(3)}}}}}\\
&\le
\begin{array}[t]{l}
\sup_{\lsubstack{\check{x}\in X,\check{y}\in Y,\check{z}\in Z\\\mathrlap{d(\check{x},\check{y})< 2\de}}}\Big(\\
\int\p{\vv{(k(h,\inv{(\cdot)}\check{y})-k(h,\inv{(\cdot)}\check{x}))\1_{K}(\cdot)}_{\infty}+ \vv{k(h,\inv{(\cdot)}\check{z})}_{\infty}\frac*{\vv{\chih-\chi}}{\vv{\chi}}}\1_F\p{h}dh\\
+\int\p{\vv{(k(h,\inv{(\cdot)}\check{y})-k(h,\inv{(\cdot)}\check{x}))\1_{K}(\cdot)}_{\infty}+ \vv{k(h,\inv{(\cdot)}\check{z})}_{\infty}\frac*{\vv{\chih-\chi}}{\vv{\chi}}}\1_{G\setminus F}\p{h}dh\Big)
\end{array}\\
&<\frac{2}{6}\ep+\frac{4}{6}\ep\:\qedhere
\end{align*}
\end{proof}
\begin{rmk}\label{rmk:po(4.1)}
Theorem~\ref{thm:po(4)} will be improved provided the group action and $k$ satisfy a Lipschitz continuity condition, $k$ satisfies a linear decay condition, and the group has strict polynomial growth.
\end{rmk}
\begin{lem}\label{lem:po(5)}
\newcommand{\chih}{\cramped{\chi^h}}%
Assume $ \ab{G,\lambda} $ is unimodular.
Let $ r\ge 1 \:$
Define
\[\ts
\chi_r=\chi:G\to\lbrp{0,+\infty}:g\mapsto(\frac{r-\v{g}}{r})\1_{B(r)}(g)\:
\]
Assume $ \ab{(G,\ell),\lambda} $ is strictly-polynomial-growing:
	there exist $ b\ge 1$ and $ c_1\ge c_0>0 $ such that for each $ r\ge 1 \;$ $ c_0r^b\le \lambda\p{B(r)}\le c_1r^b \:$
Fix $ h\in G \:$
Then
\[\ts
\vv{\chih-\chi}
\le\min\s{\frac{\v{h}}{r}(\frac{(2+b)^{2+b}c_1}{2b^bc_0})^{\frac{1}{2}},2^{\frac{1}{2}}}\vv{\chi}\;
\]
where $ \chih:G\to\lbrp{0,+\infty}: g\mapsto\chi\p{gh} \:$
\end{lem}
\begin{proof}
\newcommand{\chih}{\cramped{\chi^h}}%
Define $ r'\defeq (\frac{2}{2+b})r \:$
Observe
\begin{align*}
&\vv{\chih-\chi}^2\\
&=\int_{B(r)\inv{h}\setminus B(r)}(\frac{r-\v{gh}}{r})^2dg
+\int_{B(r)\setminus B(r)\inv{h}}(\frac{r-\v{g}}{r})^2dg
+\int_{B(r)\cap B(r)\inv{h}}(\frac{\v{g}-\v{gh}}{r})^2dg\\
&\le\int_{B(r)\inv{h}\setminus B(r)}(\frac{\v{g}-\v{gh}}{r})^2dg
+\int_{B(r)\setminus B(r)\inv{h}}(\frac{\v{gh}-\v{g}}{r})^2dg
+\int_{B(r)\cap B(r)\inv{h}}(\frac{\v{g}-\v{gh}}{r})^2dg\\
&\le\int_{B(r)\inv{h}\setminus B(r)}r^{-2}\v{h}^2dg
+\int_{B(r)\setminus B(r)\inv{h}}r^{-2}\v{h}^2dg
+\int_{B(r)\cap B(r)\inv{h}}r^{-2}\v{h}^2dg\\
&\le r^{-2}\v{h}^22(c_1r^b)
\shortintertext{and}
&\vv{\chi}^2\\
&=\int_{B(r)\setminus B(r-r')}(\frac{r-\v{g}}{r})^2dg
+\int_{B(r-r')}(\frac{r-\v{g}}{r})^2dg\\
&\ge\int_{B(r-r')}(\frac{r-\v{g}}{r})^2dg\\
&\ge\int_{B(r-r')}r^{-2}(r')^2dg\\
&\ge r^{-2}(r')^2(c_0(r-r')^b)\:
\end{align*}
Also,
\[\ts
\vv{\chih-\chi}^2\le \vv{\chih}^2+\vv{\chi}^2=2\vv{\chi}^2\:\qedhere
\]
\end{proof}
\begin{rmk}\label{rmk:po(5.1)}
\allowdisplaybreaks%
\leavevmode%
\newcommand{\chih}{\cramped{\chi^h}}%
\begin{enumerate}[(a)]
\item\label{351a}
For each nonempty compact subset $ F $ of $ G \;$
	$ \sup_{h\in F}{\frac{\vv{\chih-\chi}}{\vv{\chi}}}\to 0 $ as $ r\to+\infty \:$
\item\label{351b}
If $ G $ is discrete, then
\begin{align*}\swapabovedisplayskip
\vv{\1_{B(r)}}^2=&\#B(r)\;\\
\vv{\1_{B(r)\inv{h}\!}-\1_{B(r)}}^2=&\#\p{B(r)\inv{h}\!\sd B(r)}\:
\end{align*}
\item\label{351c}
Define $ \acm r \defeq \max\s{n\in\Z:0\le n<r} \:$
If $ G=\Z $, then
\begin{align*}
\vv{\1_{B(r)}}^2
=&2\acm r +1\;\\
\vv{\1_{B(r)-m\!}-\1_{B(r)}}^2
=&\min\s{2\v{m},2(2\acm r +1)}\;\\
\frac{\vv{\1_{B(r)-m\!}-\1_{B(r)}}}{\vv{\1_{B(r)}}}
=&O(r^{-\frac{1}{2}})\big|_{m\in\Z,r\to+\infty}\\
\shortintertext{but}
\vv{\chi}^2
=&\frac{2}{3}(\frac{(2\acm r +1)(\acm r +1)}{2r^2}+\frac{3(r-\acm r -1)}{r})\acm r  +1\;\\
\vv{\cramped{\chi^m}-\chi}^2
=&\begin{cases}
\lsubstack{\cramped{
0}
}&\lsubstack{\iph m=0}
\\
\lsubstack{\cramped{
\frac{m^2}{r}\bigp{2-\frac{m}{r}+\frac{2(r-\acm r )^2-2(r-\acm r )+1}{mr}}}
}&\lsubstack{\iph 1\le m \le \acm r }
\\
\lsubstack{\cramped{
\frac{m}{3}\bigp{(\frac{m}{r})^2-6(\frac{m}{r})+12(\frac{\acm r }{r})}
-\frac{\acm r }{3}\bigp{4(\frac{\acm r }{r})}}\\
\cramped{\hphantom{\frac{m}{3}}
+\frac{m}{3}\bigp{6(1-(\frac{\acm r }{r})^2)+6(1-\frac{\acm r }{r})(\frac{1}{r})-(\frac{1}{r})^2}
-\frac{\acm r }{3}\bigp{4(2\acm r +1)\p{1-\frac{\acm r }{r}}(\frac{1}{r})-(\frac{1}{r})^2}
}
}&\lsubstack{\iph \acm r +1\le m \le 2\acm r }
\\
\lsubstack{\cramped{
2\vv{\chi}^2}
}&\lsubstack{\iph 2\acm r +1\le m}\;
\end{cases}\\
\frac{\vv{\cramped{\chi^m}-\chi}}{\vv{\chi}}=&O(r^{-1})\big|_{m\in\Z,r\to+\infty}\:
\end{align*}
\item
Lemma~\ref{lem:po(5)} establishes that strict polynomial growth implies amenability; see Theorem~\ref{thm:spg1} and part (\ref{351a}) above.
In the case where $ G $ is discrete, one could consider indicator functions on a Følner sequence instead of tent functions---which requires a metric to define; see part (\ref{351b}) above.
An advantage of using tent functions instead of indicator functions on balls $ B(r) $ is that the limiting behavior is of order $r^{-1}$ which is sharp; see Lemma~\ref{lem:po(5)} and part (\ref{351c}) above.
\item
In the case where $ G=\Z \;$ a 1990 result by M.~Choi; G.A.~Elliott; N.~Yui \cite{MP1990_ChoiElliottYui} establishes spectral $ \frac{1}{3} $-Hölder continuity for one-dimensional discrete quasiperiodic Schrödinger operators with Lipschitz continuous potentials, e.g., the Almost Mathieu Operator, by utilizing indicator functions.
Within the same year, J.~Avron; P.H.M.v.~Mouche; B.~Simon \cite{MP1990_AvronMoucheSimon} improved the result and established spectral $ \frac{1}{2} $-Hölder continuity by utilizing tent functions.
\end{enumerate}
\end{rmk}
\begin{thm}\label{thm:po(6)}
Assume $ \ab{G,\lambda} $ is unimodular, $ A_x $ is normal for every $ x\in Z \;$ and the following.
\begin{enumerate}[(i)]
\item
$ \alpha $ satisfies a Lipschitz continuity condition:
	there exists $ c_\alpha>0 $ such that for each $ g\in G \;$
	\[\ts d(\inv{g}x,\inv{g}y)\le (c_\alpha\v{g}+1)d(x,y)\: \]
\item
$ k $ satisfies a Lipschitz continuity condition:
	there exists $ c_k\in L^1_+(G) $ such that for $ \lambda $-a.e.\ $ h\in G \;$
	\[\ts \v{k(h,x)-k(h,y)}\le c_{k}\p{h}d(x,y)\: \]
\item
$ k $ satisfies a linear decay condition:
	there exists $ c_s>0 $ such that for each $ r\ge 0 \;$
	\[\ts \sup_{x\in Z}\int \vv{k(h,\inv{(\cdot)}x)\min\s{\v{h},r}}_\infty dh\le c_s \:\]
\item
$ \ab{(G,\ell),\lambda} $ is strictly-polynomial-growing:
	there exist $ b\ge 1$ and $ c_1\ge c_0>0 $ such that for each $ r\ge 1 \;$
	\[\ts c_0r^b\le \lambda\p{B(r)}\le c_1r^b \:\]
\end{enumerate}
Fix $ X,Y\in\Ss \:$
Then for each $ r\ge 1 \;$
\[\ts
\dist_{\ssm H}\p{\spec\p{A_X},\spec\p{A_Y}}\le \vv{c_{k}}_1(c_\alpha r+1)d_{\ssm H}(X,Y)+\frac*{c_s}{r}(\frac*{(2+b)^{2+b}c_1}{2b^bc_0})^{\frac*{1}{2}}
\]
and
\[\ts
d_{\ssm H}(X,Y)\le{\min\s{\de,1}}
\lRightarrow
\dist_{\ssm H}\p{\spec\p{A_X},\spec\p{A_Y}}\le C d_{\ssm H}(X,Y)^{\frac*{1}{2}}\;
\]
where $ \de\defeq \frac*{c_s}{\vv{c_k}_1c_{\alpha}}(\frac*{(2+b)^{2+b}c_1}{2b^bc_0})^{\frac*{1}{2}} $ and $ C\defeq 2(\vv{c_k}_1c_{\alpha}c_s(\frac*{(2+b)^{2+b}c_1}{2b^bc_0})^{\frac*{1}{2}})^{\frac*{1}{2}}+\vv{c_k}_1 \:$
\end{thm}
\begin{proof}
\newcommand{\chih}{\cramped{\chi^h}}%
Fix $ r\ge 1 \:$
Define $\chi_r=\chi:G\to\lbrp{0,+\infty}:g\mapsto(\frac*{r-\v{g}}{r})\1_{B(r)}(g)\:$
Let $ \rho>0 \:$
Observe
\begin{align*}
&\dist_{\ssm H}\p{\spec\p{A_X},\spec\p{A_Y}}\\
&\le\sup_{\lsubstack{\check{x}\in X,\check{y}\in Y,\check{z}\in Z\\\mathrlap{d(\check{x},\check{y})< d_{\ssm H}(X,Y)+\rho}}}
	\int\vv{(k(h,\inv{(\cdot)}\check{y})-k(h,\inv{(\cdot)}\check{x}))\1_{\supp\p{\chi}}(\cdot)}_{\infty}+ \vv{k(h,\inv{(\cdot)}\check{z})}_{\infty}\frac*{\vv{\chih-\chi}}{\vv{\chi}}dh
\tags{\p{\lsubstack{\textrm{Lem~\ref{lem:po(3)}}}}}\\
&\le\sup_{\lsubstack{\check{x}\in X,\check{y}\in Y,\check{z}\in Z\\\mathrlap{d(\check{x},\check{y})< d_{\ssm H}(X,Y)+\rho}}}
	\int\vv{
	c_{k}\p{h}\mspace{1mu minus 1mu}d(\inv{(\cdot)}\check{y},\inv{(\cdot)}\check{x})
	\1_{\supp\p{\chi}}(\cdot)}_{\infty}+ \vv{k(h,\inv{(\cdot)}\check{z})}_{\infty}\frac*{\vv{\chih-\chi}}{\vv{\chi}}dh
\\
&\le\sup_{\lsubstack{\check{x}\in X,\check{y}\in Y,\check{z}\in Z\\\mathrlap{d(\check{x},\check{y})< d_{\ssm H}(X,Y)+\rho}}}
	\int
	c_{k}\p{h}(c_\alpha r+1)d(\check{y},\check{x})
	+ \vv{k(h,\inv{(\cdot)}\check{z})}_{\infty}\frac*{\vv{\chih-\chi}}{\vv{\chi}}dh
\\
&\le
	\vv{c_{k}}_1(c_\alpha r+1)(d_{\ssm H}(X,Y)+\rho)
	+\sup_{\lsubstack{\check{z}\in Z}}\int \vv{k(h,\inv{(\cdot)}\check{z})}_{\infty}\min\s{\frac*{\v{h}}{r}(\frac*{(2+b)^{2+b}c_1}{2b^bc_0})^{\frac*{1}{2}},2^{\frac*{1}{2}}}dh
\tags{\p{\lsubstack{\textrm{Lem~\ref{lem:po(5)}}}}}\\
&\le
	\vv{c_{k}}_1(c_\alpha r+1)(d_{\ssm H}(X,Y)+\rho)
	+\frac*{c_s}{r}(\frac*{(2+b)^{2+b}c_1}{2b^bc_0})^{\frac*{1}{2}}\:
\end{align*}
Here $ \chih:G\to \lbrp{0,+\infty}:g\mapsto \chi(gh) \:$
Also, if
\[\ts d_{\ssm H}(X,Y)\le \min\s{\de,1} \;\]
then $1\le (d_{\ssm H}(X,Y)^{-1}\de)^{\frac*{1}{2}}$ and $ d_{\ssm H}(X,Y)\le d_{\ssm H}(X,Y)^{\frac*{1}{2}} $ and
\begin{align*}
&\dist_{\ssm H}\p{\spec\p{A_X},\spec\p{A_Y}}\\
&\le \vv{c_{k}}_1(c_\alpha r+1)d_{\ssm H}(X,Y)+\frac*{c_s}{r}(\frac*{(2+b)^{2+b}c_1}{2b^bc_0})^{\frac*{1}{2}}\Big|_{r=(d_{\ssm H}(X,Y)^{-1}\de)^{\frac*{1}{2}}}\\
&\le C d_{\ssm H}(X,Y)^{\frac*{1}{2}}\:\qedhere
\end{align*}
\end{proof}
\begin{thm}\label{thm:po(7)}
Assume $ \ab{G,\lambda} $ is unimodular, $ A_x $ is normal for every $ x\in Z \;$ and the following.
\begin{enumerate}[(i)]
\item
$ \alpha $ satisfies a Lipschitz continuity condition:
	there exists $ c_\alpha>0 $ such that for each $ g\in G \;$
	\[\ts d(\inv{g}x,\inv{g}y)\le (c_\alpha\v{g}+1)d(x,y)\: \]
\item
$ k $ satisfies a locally-constant condition:
	there exists $ c_k\in L^\infty_+(G) $ such that for $ \lambda $-a.e.\ $ h\in G \;$
	\[\ts d(x,y)<\frac*{1}{c_{k}\p{h}}\lRightarrow k(h,x)=k(h,y)\: \]
\item
$ k $ satisfies a linear decay condition:
	there exists $ c_s>0 $ such that for each $ r\ge 0 \;$
	\[\ts \sup_{x\in Z}\int \vv{k(h,\inv{(\cdot)}x)\min\s{\v{h},r}}_\infty dh\le c_s \:\]
\item
$ \ab{(G,\ell),\lambda} $ is strictly-polynomial-growing:
	there exist $ b\ge 1$ and $ c_1\ge c_0>0 $ such that for each $ r\ge 1 \;$
	\[\ts c_0r^b\le \lambda\p{B(r)}\le c_1r^b \:\]
\end{enumerate}
Fix $ X,Y\in\Ss \:$
Then
\[\ts
d_{\ssm H}(X,Y)\le \de\lRightarrow
\dist_{\ssm H}\p{\spec\p{A_X},\spec\p{A_Y}}\le C d_{\ssm H}(X,Y)\;
\]
where $ \de\defeq \frac*{1}{(\vv{c_k}_\infty+1)(c_\alpha+1)}$ and $ C\defeq (\vv{c_k}_\infty+1)(c_\alpha+1)c_s(\frac*{(2+b)^{2+b}c_1}{2b^bc_0})^{\frac*{1}{2}} \:$
\end{thm}
\begin{proof}
\newcommand{\chih}{\cramped{\chi^h}}%
Assume $ d_{\ssm H}(X,Y)\le \de \:$
Observe
\begin{align*}
\frac*{1}{(\vv{c_k}_\infty+1)\de c_\alpha}-\frac*{1}{c_\alpha}=&1\;\\
\frac*{c_\alpha}{1-(\vv{c_k}_\infty+1)d_{\ssm H}(X,Y)}\le& c_\alpha +1\:
\end{align*}
Define $ r\defeq \frac*{1}{(\vv{c_k}_\infty+1)d_{\ssm H}(X,Y) c_\alpha}-\frac*{1}{c_\alpha} \:$
Observe
\begin{align*}
r\ge& 1\;\\
\frac{1}{r}\le& (\vv{c_k}_{\infty}+1)(c_\alpha+1)d_{\ssm H}(X,Y)\:
\end{align*}
Define, for $ h\in G \;$
\begin{align*}\swapabovedisplayskip
\chi_r=\chi
:&G\to\lbrp{0,+\infty}:g\mapsto(\frac*{r-\v{g}}{r})\1_{B(r)}(g)\;\\
\chih
:&G\to\lbrp{0,+\infty}:g\mapsto\chi\p{gh}\:
\end{align*}
Define $ \rho\defeq \frac*{d_{\ssm H}(X,Y)}{\vv{c_k}_\infty} \:$
Observe
\[\ts
d(x,y)<d_{\ssm H}(X,Y)+\rho
\lRightarrow
\vv{d(\inv{(\cdot)}y,\inv{(\cdot)}x)\1_{\supp\p{\chi}}(\cdot)}_{\infty}
\le(c_\alpha r+1)d(x,y)
<\frac*{1}{\vv{c_k}_\infty}
\]
and
\begin{align*}\swapabovedisplayskip
&\dist_{\ssm H}\p{\spec\p{A_X},\spec\p{A_Y}}\\
&\le\sup_{\lsubstack{\check{x}\in X,\check{y}\in Y,\check{z}\in Z\\\mathrlap{d(\check{x},\check{y})< d_{\ssm H}(X,Y)+\rho}}}
	\int\vv{(k(h,\inv{(\cdot)}\check{y})-k(h,\inv{(\cdot)}\check{x}))\1_{\supp\p{\chi}}(\cdot)}_{\infty}+ \vv{k(h,\inv{(\cdot)}\check{z})}_{\infty}\frac*{\vv{\chih-\chi}}{\vv{\chi}}dh
\tags{\p{\lsubstack{\textrm{Lem~\ref{lem:po(3)}}}}}\\
&\le
	0+ \sup_{\lsubstack{\check{z}\in Z}}\int\vv{k(h,\inv{(\cdot)}\check{z})}_{\infty}\frac*{\vv{\chih-\chi}}{\vv{\chi}}dh
	\\
&\le
\sup_{\lsubstack{\check{z}\in Z}}\int \vv{k(h,\inv{(\cdot)}\check{z})}_{\infty}\min\s{\frac*{\v{h}}{r}(\frac*{(2+b)^{2+b}c_1}{2b^bc_0})^{\frac*{1}{2}},2^{\frac*{1}{2}}}dh
\tags{\p{\lsubstack{\textrm{Lem~\ref{lem:po(5)}}}}}\\
&\le
\frac*{c_s}{r}(\frac*{(2+b)^{2+b}c_1}{2b^bc_0})^{\frac*{1}{2}}
\\
&\le C d_{\ssm H}(X,Y)\:\qedhere
\end{align*}
\end{proof}
	\clearpage
\section{Proof of Lemma \ref{lem:d-dof1}}\label{sec:A1}
{
\noindent
Let $ \ab{G,\lambda} $ be a second-countable locally-compact Hausdorff group endowed with a left Haar measure.
Let $ \ab{Z,(G,\alpha)} $ be a topological space endowed with a left continuous action.
Let $ k:G\times Z\to\C \;$ where $ k $ is measurable and
	$ \sup_{x\in Z}\int \vv{k(h,\inv{(\cdot)}x)}_{\infty} dh<+\infty\: $
Define, for $ x\in Z \;$
	${A}_x:L^2(G)\to L^2(G)$ by
\[\ts
(A_x\psi)(g)=\int k(h,\inv{g}x)\psi(gh) dh=\int k(\inv{g}h,\inv{g}x)\psi(h)dh
\]
for every $ \psi\in L^2(G)\, g\in G \:$
Let $ \Ss $ be the collection of all subsystems of $ \ab{Z,(G,\alpha)} \:$
Define, for $ X\in\Ss \;$
	$ \spec(A_X)\defeq \overline{\bigcup_{x\in X}\spec(A_x)} \:$
}
\begin{prop}[Lemma \ref{lem:d-dof1}(\ref{d-dof1b})]\label{prop:A1(1)}
Assume $ \ab{G,\lambda} $ is unimodular.
Fix $ x,y\in Z $ and $ \psi\in L^2(G) $ and $ j\in G \:$
Then
\begin{align*}\swapabovedisplayskip
\vv{A_x\psi}\le&
	(\int\vv{k(h,\inv{(\cdot)}x)\1_{\supp(\psi)\inv{h}}(\cdot)}_{\infty}dh)\vv{\psi}\;\\
\vv{A_x\!{}^{*}\psi}\le&
	(\int\vv{k(h,\inv{(\cdot)}x)\1_{\supp(\psi)}(\cdot)}_{\infty}dh)\vv{\psi}\;\\
\vv{(A_y-A_x)\psi}\le&
	(\int\vv{\p{k(h,\inv{(\cdot)}y)-k(h,\inv{(\cdot)}x)}\1_{\supp(\psi)\inv{h}}(\cdot)}_{\infty}dh)\vv{\psi}\;\\
\vv{(A_y-A_x)^{*}\psi}\le&
	(\int\vv{\p{k(h,\inv{(\cdot)}y)-k(h,\inv{(\cdot)}x)}\1_{\supp(\psi)}(\cdot)}_{\infty}dh)\vv{\psi}\;\\
U_{\inv{j}}A_xU_{j}=&A_{\inv{j}x}\:
\end{align*}
\end{prop}
\begin{proof}
Observe
\begin{align*}
&\vv{A_x\psi}^2\\
&=\int\v{\int k(h,\inv{g}x)\psi(gh) dh}^2dg\\
&\le\int\p{\int \v{k(h,\inv{g}x)}\v{\psi(gh)} dh}^2dg\\
&\le\int
\p{\int \v{k(h,\inv{g}x)}\1_{\supp\p{\psi}}(gh) dh}
\p{\int \v{k(h,\inv{g}x)}\v{\psi(gh)}^2 dh}dg
\tagps{\textrm{Hölder's}\\\textrm{inequality}}\\
&\le\int
\p{\int \vv{k(h,\inv{(\cdot)}x)\1_{\supp(\psi)\inv{h}}(\cdot)}_{\infty} dh}
\p{\int \vv{k(h,\inv{(\cdot)}x)\1_{\supp(\psi)\inv{h}}(\cdot)}_{\infty}\v{\psi(gh)}^2 dh}dg\\
&= \p{\int \vv{k(h,\inv{(\cdot)}x)\1_{\supp(\psi)\inv{h}}(\cdot)}_{\infty} dh}^2\vv{\psi}^2
\shortintertext{and}
&\vv{A_x\!{}^{*}\psi}^2\\
&=\int\v{\int \overline{k({\inv{h}},{\inv{h}\inv{g}}x)}\psi(gh) dh}^2dg\\
&\le\int\p{\int \v{\overline{k({\inv{h}},{\inv{h}\inv{g}}x)}}\v{\psi(gh)} dh}^2dg\\
&\le\int
\p{\int \v{k({\inv{h}},{\inv{(gh)}}x)}\1_{\supp\p{\psi}}(gh) dh}
\p{\int \v{k({\inv{h}},{\inv{(gh)}}x)}\v{\psi(gh)}^2 dh}dg\\
&\le\int
\p{\int \vv{k({\inv{h}},{\inv{(\cdot)}}x)\1_{\supp\p{\psi}}(\cdot)}_{\infty} dh}
\p{\int \vv{k({\inv{h}},{\inv{(\cdot)}}x)\1_{\supp\p{\psi}}(\cdot)}_{\infty}\v{\psi(gh)}^2 dh}dg\\
&= \p{\int \vv{k({\inv{h}},{\inv{(\cdot)}}x)\1_{\supp\p{\psi}}(\cdot)}_{\infty} dh}^2\vv{\psi}^2\\
\shortintertext{and}
&\vv{(A_y-A_x)\psi}^2\\
&=\int\v{\int\p{k(h,{\inv{g}}y)-k(h,{\inv{g}}x)}\psi(gh) dh}^2dg\\
&\le\p{\int \vv{\p{k(h,{\inv{(\cdot)}}y)-k(h,{\inv{(\cdot)}}x)}\1_{\supp(\psi)\inv{h}}(\cdot)}_{\infty} dh}^2\vv{\psi}^2
\shortintertext{and}
&\vv{(A_y-A_x)^*\psi}^2\\
&=\int\v{\int \p{\overline{k({\inv{h}},{\inv{h}\inv{g}}y)}-\overline{k({\inv{h}},{\inv{h}\inv{g}}x)}}\psi(gh) dh}^2dg\\
&\le \p{\int \vv{\p{k({\inv{h}},\inv{(\cdot)}y)-k({\inv{h}},\inv{(\cdot)}x)}\1_{\supp\p{\psi}}(\cdot)}_{\infty} dh}^2\vv{\psi}^2\:
\end{align*}
Also,
\[\ts
U_{\inv{j}}A_xU_{j}
=\int d\lambda(h)U_{\inv{j}}[k(h,\inv{(\cdot)}x)]U_{j}R_h
=\int d\lambda(h)[k(h,\inv{(\cdot)}\inv{j}x)]R_h
=A_{\inv{j}x}\:\qedhere
\]
\end{proof}
\begin{prop}[Lemma \ref{lem:d-dof1}(\ref{d-dof1c})]\label{prop:A1(2)}
Assume $ \ab{G,\lambda} $ is unimodular, $ A_x $ is normal for every $ x\in Z \;$ and the following.
\begin{enumerate}[(i)]
\item\label{A1(2)1}
$ k $ satisfies a continuity condition:
	for each $ x\in Z $ and for each $ \ep>0 $ and for each nonempty compact subset $ F $ of $ G \;$
	$ {\int\v{k(h,y)-k(h,x)}\1_F\p{h}dh}<\ep $ for every $ y $ in some neighborhood of $ x \:$
\item\label{A1(2)2}
\begin{iws}{-0.059pt}
$ k $ satisfies a decay condition:
	for each $ x\in Z $ and for each $ \ep>0 \;$ there exists a~nonempty compact subset $ F $ of $ G $ such that
	$ \int\vv{k(h,\inv{(\cdot)}y)\1_{G\setminus F}(h)}_{\infty}dh<\ep $ for every $ y $ in some neighborhood of $ x \:$
\end{iws}
\end{enumerate}
Fix $ x\in Z \:$
Then for each $ \ep>0 \;$ $ \sup_{E\in\spec(A_x)}\dist(E,\spec(A_y))<\ep $ for every $ y $ in some neighborhood of $ x \;$ i.e.,
\[\ts \lim_{y\xrightarrow[]{\smash{Z}}x}\sup_{E\in\spec({A}_x)}\dist(E,\spec({A}_y))=0\: \]
\end{prop}
\begin{proof}
It suffices to show that
	for each $ E\in \spec\p{A_x} $ and for each $ \ep>0 \;$ $ \dist(E,\spec(A_y))<\ep $ for every $ y $ in some neighborhood of $ x $
	since $ \spec\p{A_x} $ is compact.
Fix $ E\in\spec\p{A_x}$ and $ \ep>0 \:$
Without loss of generality, assume $ E\notin \spec\p{A_y} \:$
Observe there exists a nonempty compact subset $ F $ of $ G $ such that
\[\ts
\int\vv{k(h,\inv{(\cdot)}y)\1_{G\setminus F}(h)}_{\infty}dh<\frac{1}{4}\ep
\]
for every $ y $ in some neighborhood $ U $ of $ x $ since $ k $ satisfies a decay condition.
Observe $ A_x $ is normal and, by Weyl's criterion, there exists a nonzero vector $ \vp $ such that $ \vv{(A_x -E)\vp}< \frac{1}{4}\ep\vv{\vp} \:$
Without loss of generality, assume $ \vp $ is both continuous and compactly-supported.
Define $ K\defeq \supp\p{\vp} \:$
Observe
\[\ts
\int\vv{\p{k(h,\inv{(\cdot)}y)-k(h,\inv{(\cdot)}x)}\1_K(\cdot)}_\infty\1_F\p{h}dh<\frac{1}{4}\ep
\]
for every $ y $ in some neighborhood $ U' $ of $ x $
since the group action is continuous and $ k $ satisfies a continuity condition.
Observe for each $ y\in U \cap U' \;$ $ A_y $ is normal and
\[\ts
\vv{(A_y-E)^{-1}}^{-1}\vv{\vp}
\le
\vv{(A_y-E)\vp}
=\vv{(A_y-E)^*\vp}
< (\vv{(A_y-A_x)^*\frac{\vp}{\vv{\vp}}}+\frac{1}{4}\ep)\vv{\vp}
\]
and
\begin{align*}\swapabovedisplayskip
&\dist\p{E,\spec\p{A_y}}\\
&=\vv{(A_y-E)^{-1}}^{-1}\\
&<\vv{(A_y-A_x)^*\frac{\vp}{\vv{\vp}}}+\frac{1}{4}\ep\\
&\le\int\vv{(k(h,\inv{(\cdot)}y)-k(h,\inv{(\cdot)}x))\1_{K}(\cdot)}_{\infty}dh+\frac{1}{4}\ep
\tagps{\textrm{Prop~\ref{prop:A1(1)}}}\\
&<\int\vv{(k(h,\inv{(\cdot)}y)-k(h,\inv{(\cdot)}x))\1_{K}(\cdot)}_{\infty}\1_{F}(h)dh+\frac{3}{4}\ep\\
&< \ep\:\qedhere
\end{align*}
\end{proof}
\begin{rmk}\label{rmk:A1(2.1)}
\leavevmode%
\begin{enumerate}[(a)]
\item
To see how conditions \textit{(\ref{A1(2)1})} and \textit{(\ref{A1(2)2})} follow from the setting for Lemma~\ref{lem:d-dof1},
	recall the setting has $ f_{(\cdot)}\in C_c\p{G,C_b(Z)} $ and $ f_h(x)=k(h,x) \:$
\item
To see how $ \spec\p{A_x}=\spec\p{A_{X\p{x}}} $ follows from Proposition~\ref{prop:A1(1)} and Proposition~\ref{prop:A1(2)},
	if $ y =\inv{j}x $ for some $ j\in G \;$ then $ \spec\p{A_{y}}=\spec\p{A_x} $ since
		$ A_{y}=U_{\inv{j}}A_xU_j $ and
	if $ y $ is a limit-point of the orbit of $ x \;$ then $ \spec\p{A_y}\subseteq \spec\p{A_x} $ since
	$ \inf_{j\in G}\sup_{E\in\spec\p{A_y}}\dist\p{E,\spec\p{A_{\inv{j}x}}}=0 \:$
\end{enumerate}
\end{rmk}
\clearpage
\bookmarksetup{startatroot}
\subsection*{Acknowledgments}
We thank Anton Gorodetski for his guidance and support.
We thank Jean Bellissard and Horia Cornean for discussions on Lipschitz partitions of unity and amenability.
We thank Svetlana Jitomirskaya and David Damanik for commenting on how to improve the initial draft of this paper.
We thank Yuji Tachikawa, Yasuyuki Hatsuda, and Hosho Katsura for discussions on Hofstadter's butterfly and for providing Figure~\ref{fig:1} (see \cite{P2016_HatsKatsTach}) together with additional curves.
We thank Wencai Liu for hosting the Spectral Theory Reading Seminar (\href{https://sites.google.com/view/reading-seminar/home}{2021}).
\subsection*{Acknowledgment of Support}
This project was supported by the Deutsche Forschungsgemeinschaft [\href{https://gepris.dfg.de/gepris/projekt/412141125?language=en}{BE 6789/1-1 to S.B.}] and the second author was partially supported by NSF grant \href{https://www.nsf.gov/awardsearch/showAward?AWD_ID=1855541&HistoricalAwards=false}{DMS-1855541} (PI - A.~Gorodetski) and NSF fellowship award \href{https://www.nsf.gov/awardsearch/showAward?AWD_ID=2213277&HistoricalAwards=false}{DMS-2213277}.
\subsection*{Transparency-and-Accountability Declarations}
The authors have no relevant financial or non-financial interests to disclose beyond the above acknowledgment of support and
the authors have no competing interests to declare that are relevant to the content of this article.
\NEWvbadness{10000}
\microtypesetup{deactivate}


\begin{thebibliography}{00}
\bibitem{P1964_Azbel}
	M.~Azbel, \textit{Energy spectrum of a conduction electron in a magnetic field}, Sov.\ Phys.\ JETP \textbf{19}, (1964), 634--645.
\bibitem{MP1981_ASLim}
	J.~Avron and B.~Simon, \textit{Almost periodic Schrödinger operators.\ I.\ Limit periodic potentials}, Comm.\ Math.\ Phys.\ \textbf{82}, no.\ 1, (\href{https://mathscinet.ams.org/mathscinet-getitem?mr=638515}{1981}), 101--120.
\bibitem{MP1983_AvronSimonAMO}
	J.~Avron and B.~Simon, \textit{Almost periodic Schrödinger operators.\ II.\ The integrated density of states}, Duke Math.\ J.\ \textbf{50}, no.\ 1, (\href{https://mathscinet.ams.org/mathscinet-getitem?mr=700145}{1983}), 369--391.
\bibitem{MP1985_AvronSimonGCont}
	J.~Avron and B.~Simon, \textit{Stability of gaps for periodic potentials under variation of a magnetic field}, J.\ Phys.\ A \textbf{18}, no.\ 12, (\href{https://mathscinet.ams.org/mathscinet-getitem?mr=804317}{1985}), 2199--2205.
\bibitem{P1990_AustinWilkinson}
	E.J.~Austin and M.~Wilkinson, \textit{Phase space lattices with threefold symmetry}, J.\ Phys.\ A \textbf{23}, no.\ 12, (1990), 2529--2554.
\bibitem{MP1990_AvronMoucheSimon}
	J.~Avron and P.H.M.v~Mouche and B.~Simon, \textit{On the measure of the spectrum for the Almost Mathieu Operator}, Comm.\ Math.\ Phys.\ \textbf{132}, no.\ 1, (\href{https://mathscinet.ams.org/mathscinet-getitem?mr=1069202}{1990}), 103--118.
\bibitem{P1994_AustinWilkinson}
	E.J.~Austin and M.~Wilkinson, \textit{Spectral dimension and dynamics for Harper's equation}, Phys.\ Rev.\ B \textbf{50}, no.\ 3, (1994), 1420--1429.
\bibitem{MP2009_ALim}
	A.~Avila, \textit{On the spectrum and Lyapunov exponent of limit periodic Schrödinger operators}, Comm.\ Math.\ Phys.\ \textbf{288}, no.\ 3, (\href{https://mathscinet.ams.org/mathscinet-getitem?mr=2504859}{2009}), 907--918.
\bibitem{MP2009_AvilaJitomirskayaAMO}
	A.~Avila and S.~Jitomirskaya, \textit{The Ten Martini Problem}, Ann.\ of Math.\ \textbf{170}, no.\ 1, (\href{https://mathscinet.ams.org/mathscinet-getitem?mr=2521117}{2009}), 303--342.
\bibitem{MP2020_AvniBreuerSimonLGO}
	N.~Avni and J.~Breuer and B.~Simon, \textit{Periodic Jacobi matrices on trees}, Adv.\ Math.\ \textbf{370}, (\href{https://mathscinet.ams.org/mathscinet-getitem?mr=4103777}{2020}), 42 pp.
\bibitem{MP1985_BellissardLimaTestard}
	J.~Bellissard and R.~Lima and D.~Testard, \textit{Almost periodic Schrödinger operators}, Mathematics + physics, vol.\ 1, World Sci.\ Publishing, Singapore, (\href{https://mathscinet.ams.org/mathscinet-getitem?mr=849342}{1985}), 1--64.
\bibitem{MP1986_Bellissard}
	J.~Bellissard, \textit{\(K\)-theory of \(C^*\)-algebras in solid state physics}, Statistical mechanics and field theory:\ mathematical aspects (Groningen, 1985), Lecture Notes in Phys.\ \textbf{257}, (\href{https://mathscinet.ams.org/mathscinet-getitem?mr=862832}{1986}), 99--156.
\bibitem{P1987_BellStinchcombe}
	S.C.~Bell and R.B.~Stinchcombe, \textit{Hierarchical band clustering and fractal spectra in incommensurate systems}, J.\ Phys.\ A \textbf{20}, no.\ 11, (1987), L739--L744.
\bibitem{P1990_BellissardRammal}
	J.~Bellissard and R.~Rammal, \textit{An algebraic semi-classical approach to Bloch electrons in a magnetic field}, J.\ Phys.\ France \textbf{51}, (1990), 1803--1830.
\bibitem{MP1991_BellissardIochumTestardGCont}
	J.~Bellissard and B.~Iochum and D.~Testard, \textit{Continuity properties of the electronic spectrum of 1D quasicrystals}, Comm.\ Math.\ Phys.\ \textbf{141}, no.\ 2, (\href{https://mathscinet.ams.org/mathscinet-getitem?mr=1133271}{1991}), 353--380.
\bibitem{MP1992_BellissardBovierGhezGLT}
	J.~Bellissard and A.~Bovier and J.-M.~Ghez, \textit{Gap labelling theorems for one-dimensional discrete Schrödinger operators}, Rev.\ Math.\ Phys.\ \textbf{4}, no.\ 1, (\href{https://mathscinet.ams.org/mathscinet-getitem?mr=1160136}{1992}), 1--37.
\bibitem{MP1992_BellissardGLT}
	J.~Bellissard, \textit{Gap labelling theorems for Schrödinger operators}, From number theory to physics (Les Houches, 1989), (\href{https://mathscinet.ams.org/mathscinet-getitem?mr=1221111}{1992}), 538--630.
\bibitem{MP1994_BellissardGCont}
	J.~Bellissard, \textit{Lipshitz continuity of gap boundaries for Hofstadter-like spectra}, Comm.\ Math.\ Phys.\ \textbf{160}, no.\ 3, (\href{https://mathscinet.ams.org/mathscinet-getitem?mr=1266066}{1994}), 599--613.
\bibitem{MP2001_BGSSkew}
	J.~Bourgain and M.~Goldstein and W.~Schlag, \textit{Anderson localization for Schrödinger operators on \(\Z\) with potentials given by the skew-shift}, Comm.\ Math.\ Phys.\ \textbf{220}, no.\ 3, (\href{https://mathscinet.ams.org/mathscinet-getitem?mr=1843776}{2001}), 583--621.
\bibitem{MP2002_BourgainSquareAMO1of2}
	J.~Bourgain, \textit{On the spectrum of lattice Schrödinger operators with deterministic potential}, J.\ Anal.\ Math.\ \textbf{87}, (\href{https://mathscinet.ams.org/mathscinet-getitem?mr=1945277}{2002}), 37--75.
\bibitem{MP2005_BourgainReview}
	J.~Bourgain, \textit{Green's function estimates for lattice Schrödinger operators and applications}, Annals of Mathematics Studies, \textbf{158}, Princeton University Press, Princeton, NJ, (\href{https://mathscinet.ams.org/mathscinet-getitem?mr=2100420}{2005}), x+173 pp.
\bibitem{MP2018_BeckusBellissardDeNittisGCT1of2}
	S.~Beckus and J.~Bellissard and G.~De~Nittis, \textit{Spectral continuity for aperiodic quantum systems I.\ General theory}, J.\ Funct.\ Anal.\ \textbf{275}, no.\ 11, (\href{https://mathscinet.ams.org/mathscinet-getitem?mr=3861728}{2018}), 2917--2977.
\bibitem{MP2019_BeckusBellissardCorneanGCont}
	S.~Beckus and J.~Bellissard and H.~Cornean, \textit{Hölder continuity of the spectra for aperiodic Hamiltonians}, Ann.\ Henri Poincar\'{e} \textbf{20}, no.\ 11, (\href{https://mathscinet.ams.org/mathscinet-getitem?mr=4019198}{2019}), 3603--3631.
\bibitem{MP2022_BrunoCalziLGO}
	T.~Bruno and M.~Calzi, \textit{Schrödinger operators on Lie groups with purely discrete spectrum}, Adv.\ Math.\ \textbf{404}, (\href{https://mathscinet.ams.org/mathscinet-getitem?mr=4418885}{2022}), 45pp.
\bibitem{MP1987_CyconFroeseKirschSimonReview}
	H.L.~Cycon and R.G.~Froese and W.~Kirsch and B.~Simon, \textit{Schrödinger operators with application to quantum mechanics and global geometry}, Texts and Monographs in Physics, (\href{https://mathscinet.ams.org/mathscinet-getitem?mr=883643}{1987}), x+319 pp.
\bibitem{MP1990_CarmonaLacroixReview}
	R.~Carmona and J.~Lacroix, \textit{Spectral theory of random Schrödinger operators}, Probability and its Applications, (\href{https://mathscinet.ams.org/mathscinet-getitem?mr=1102675}{1990}), xxvi+587 pp.
\bibitem{MP1990_ChoiElliottYui}
	M.~Choi and G.A.~Elliott and N.~Yui, \textit{Gauss polynomials and the rotation algebra}, Invent.\ Math.\ \textbf{99}, no.\ 2, (\href{https://mathscinet.ams.org/mathscinet-getitem?mr=1031901}{1990}), 225--246.
\bibitem{MR2962855}
	H.~Cornean and R.~Purice, \textit{On the regularity of the Hausdorff distance between spectra of perturbed magnetic Hamiltonians}, Oper.\ Theory Adv.\ Appl.\ \textbf{224}, (\href{https://mathscinet.ams.org/mathscinet-getitem?mr=2962855}{2012}), 55--66.
\bibitem{MP2016_DamanikGorodetskiYessFH}
	D.~Damanik and A.~Gorodetski and W.~Yessen,\textit{The Fibonacci Hamiltonian}, Invent.\ Math.\ \textbf{206}, no.\ 3, (\href{https://mathscinet.ams.org/mathscinet-getitem?mr=3573970}{2016}), 629--692.
\bibitem{MP1982_ElliottAMO}
	G.A.~Elliott, \textit{Gaps in the spectrum of an almost periodic Schrödinger operator}, C.\ R.\ Math.\ Rep.\ Acad.\ Sci.\ Canada \textbf{4}, no.\ 5, (\href{https://mathscinet.ams.org/mathscinet-getitem?mr=675127}{1982}), 255--259.
\bibitem{MP1992_PasturFigotinReview}
	A.~Figotin and L.~Pastur, \textit{Spectra of random and almost-periodic operators}, Grundlehren der mathematischen Wissenschaften [Fundamental Principles of Mathematical Sciences] \textbf{297}, (\href{https://mathscinet.ams.org/mathscinet-getitem?mr=1223779}{1992}), viii+587 pp.
\bibitem{GR2016_Folland}
	G.B.~Folland, \textit{A course in abstract harmonic analysis}, Second edition, Textbooks in Mathematics, CRC Press, Boca Raton, FL, (\href{https://mathscinet.ams.org/mathscinet-getitem?mr=3444405}{2016}), xiii+305 pp.
\bibitem{P1991_GeiselKetzmerickPetschel}
	T.~Geisel and R.~Ketzmerick and G.~Petschel, \textit{New class of level statistics in quantum systems with unbounded diffusion}, Phys.\ Rev.\ Lett.\ \textbf{66}, no.\ 13, (1991), 1651--1654.
\bibitem{P1998_GeiselKetzmerickKruseSteinbach}
	T.~Geisel and R.~Ketzmerick and K.~Kruse and F.~Steinbach, \textit{Covering property of Hofstadter's butterfly}, Phys.\ Rev.\ B \textbf{58}, no.\ 15, (1998), 9881--9885.
\bibitem{MP2010_GLim}
	Z.~Gan, \textit{An exposition of the connection between limit-periodic potentials and profinite groups}, Math.\ Model.\ Nat.\ Phenom.\ \textbf{5}, no.\ 4, (\href{https://mathscinet.ams.org/mathscinet-getitem?mr=2662454}{2010}), 158--174.
\bibitem{MR4366010}
	M.~Gerhold and O.M.~Shalit, \textit{Dilations of q-commuting unitaries}, Int.\ Math.\ Res.\ Not.\ IMRN, no.\ 1, (\href{https://mathscinet.ams.org/mathscinet-getitem?mr=4366010}{2022}), 63--88.
\bibitem{P1955_Harper1of2}
	P.G.~Harper, \textit{Single Band Motion of Conduction Electrons in a Uniform Magnetic Field}, Proc.\ Phys.\ Soc.\ A \textbf{68}, (1955), 874--878.
\bibitem{P1955_Harper2of2}
	P.G.~Harper, \textit{The General Motion of Conduction Electrons in a Uniform Magnetic Field, with Application to the Diamagnetism of Metals}, Proc.\ Phys.\ Soc.\ A \textbf{68}, (1955), 879--892.
\bibitem{P1976_Hofstadter}
	D.~Hofstadter, \textit{Energy levels and wave functions of Bloch electrons in rational and irrational magnetic fields}, Phys.\ Rev.\ B \textbf{14}, (1976), 2239--2249.
\bibitem{MP1990_HelfferSjostrandAMO}
	B.~Helffer and J.~Sjöstrand, \textit{Analyse semi-classique pour l'équation de Harper.\ II.\ Comportement semi-classique près d'un rationnel}, M\'em.\ Soc.\ Math.\ France, no.\ 40, (\href{https://mathscinet.ams.org/mathscinet-getitem?mr=1052373}{1990}), 1--139.
\bibitem{P2016_HatsKatsTach}
	Y.~Hatsuda and H.~Katsura and Y.~Tachikawa, \textit{Hofstadter's butterfly in quantum geometry}, New J.\ Phys.\ \textbf{18}, (2016), 103023.
\bibitem{MP2019_HelfferLiuQuZhouAMO}
	B.~Helffer and Q.~Liu and Y.~Qu and Q.~Zhou, \textit{Positive Hausdorff dimensional spectrum for the critical Almost Mathieu Operator}, Comm.\ Math.\ Phys.\ \textbf{368}, no.\ 1, (\href{https://mathscinet.ams.org/mathscinet-getitem?mr=3946411}{2019}), 369--382.
\bibitem{MP1998_JitomirskayaLastAMO}
	S.~Jitomirskaya and Y.~Last, \textit{Anderson localization for the almost Mathieu equation.\ III.\ Semi-uniform localization, continuity of gaps, and measure of the spectrum}, Comm.\ Math.\ Phys.\ \textbf{195}, no.\ 1, (\href{https://mathscinet.ams.org/mathscinet-getitem?mr=1637389}{1998}), 1--14.
\bibitem{MP2002_JitomirskayaKrasAMO}
	S.~Jitomirskaya and I.~Krasovsky, \textit{Continuity of the measure of the spectrum for discrete quasiperiodic operators}, Math.\ Res.\ Lett.\ \textbf{9}, no.\ 4, (\href{https://mathscinet.ams.org/mathscinet-getitem?mr=1928861}{2002}), 413--421.
\bibitem{MP2014_JitomirskayaMaviAMO}
	S.~Jitomirskaya and R.~Mavi, \textit{Continuity of the measure of the spectrum for quasiperiodic Schrödinger operators with rough potentials}, Comm.\ Math.\ Phys.\ \textbf{325}, no.\ 2, (\href{https://mathscinet.ams.org/mathscinet-getitem?mr=3148097}{2014}), 585--601.
\bibitem{MP2019_JitomirskayaKrasAMO}
	S.~Jitomirskaya and I.~Krasovsky, \textit{Critical almost Mathieu operator:\ hidden singularity, gap continuity, and the Hausdorff dimension of the spectrum}, Ann.\ of Math.\ (\href{https://annals.math.princeton.edu/articles/18188}{to appear}).
\bibitem{MP2022_JitomirskayaZhangAMO}
	S.~Jitomirskaya and S.~Zhang, \textit{Quantitative continuity of singular continuous spectral measures and arithmetic criteria for quasiperiodic Schrödinger operators}, J.\ Eur.\ Math.\ Soc.\ \textbf{24}, no.\ 5, (\href{https://mathscinet.ams.org/mathscinet-getitem?mr=4404788}{2022}), 1723--1767.
\bibitem{P1986_KohmotoTang}
	M.~Kohmoto and C.~Tang, \textit{Global scaling properties of the spectrum for a quasiperiodic Schrödinger equation}, Phys.\ Rev.\ B \textbf{34}, no.\ 3, (1986), 2041--2044.
\bibitem{MP1989_KirschReview}
	W.~Kirsch, \textit{Random Schrödinger operators.\ A course}, Lecture Notes in Phys.\ \textbf{345}, (\href{https://mathscinet.ams.org/mathscinet-getitem?mr=1037323}{1989}), 264--370.
\bibitem{MP2011_KrugerSkew}
	H.~Krüger, \textit{The spectrum of skew-shift Schrödinger operators contains intervals}, J.\ Funct.\ Anal.\ \textbf{262}, no.\ 3, (\href{https://mathscinet.ams.org/mathscinet-getitem?mr=2863848}{2012}), 773--810.
\bibitem{MP1994_LastAMO}
	Y.~Last, \textit{Zero measure spectrum for the Almost Mathieu Operator}, Comm.\ Math.\ Phys.\ \textbf{164}, no.\ 2, (\href{https://mathscinet.ams.org/mathscinet-getitem?mr=1289331}{1994}), 421--432.
\bibitem{GR2001_Lindenstrauss}
	E.~Lindenstrauss, \textit{Pointwise theorems for amenable groups}, Invent.\ Math.\ \textbf{146}, no.\ 2, (\href{https://mathscinet.ams.org/mathscinet-getitem?mr=1865397}{2001}), 259--295.
\bibitem{MP2007_LenzPeyerVeGAmen}
	D.~Lenz and N.~Peyerimhoff and I.~Veseli\'{c}, \textit{Groupoids, von Neumann algebras and the integrated density of states}, Math.\ Phys.\ Anal.\ Geom.\ \textbf{10}, no.\ 1, (\href{https://mathscinet.ams.org/mathscinet-getitem?mr=2340531}{2007}), 1--41.
\bibitem{MP2016_LastShamisAMO}
	Y.~Last and M.~Shamis, \textit{Zero Hausdorff dimension spectrum for the Almost Mathieu Operator}, Comm.\ Math.\ Phys.\ \textbf{348}, no.\ 3, (\href{https://mathscinet.ams.org/mathscinet-getitem?mr=3555352}{2016}), 729--750.
\bibitem{GR1955_MontgomeryZippin}
	D.~Montgomery and L.~Zippin, \textit{Topological transformation groups}, Interscience Publishers, New York-London, (\href{https://mathscinet.ams.org/mathscinet-getitem?mr=73104}{1955}), xi+282 pp.
\bibitem{GR2006_Nevo}
	A.~Nevo, \textit{Pointwise ergodic theorems for actions of groups}, Handbook of dynamical systems, vol.\ 1B, (\href{https://mathscinet.ams.org/mathscinet-getitem?mr=2186253}{2006}), 871--982.
\bibitem{GR1984_Pier}
	J.-P.~Pier, \textit{Amenable locally compact groups}, Pure and Applied Mathematics (New York), A Wiley-Interscience Publication, John Wiley \& Sons, Inc., New York, (\href{https://mathscinet.ams.org/mathscinet-getitem?mr=767264}{1984}), x+418 pp.
\bibitem{GR1988_Paterson}
	A.L.T.~Paterson, \textit{Amenability}, Mathematical Surveys and Monographs, 29, American Mathematical Society, Providence, RI, (\href{https://mathscinet.ams.org/mathscinet-getitem?mr=961261}{1988}), xx+452 pp.
\bibitem{LSA2018_Pedersen}
	G.K.~Pedersen, \textit{\(C^*\)-algebras and their automorphism groups}, Pure and Applied Mathematics (Amsterdam), Academic Press, London, (\href{https://mathscinet.ams.org/mathscinet-getitem?mr=3839621}{2018}), xviii+520 pp.
\bibitem{GR1974_Struble}
	R.A.~Struble, \textit{Metrics in locally compact groups}, Compositio Math \textbf{28}, (\href{https://mathscinet.ams.org/mathscinet-getitem?mr=348037}{1974}), 217--222.
\bibitem{MP1982_SimonRP}
	B.~Simon, \textit{Almost periodic Schrödinger operators:\ a review}, Adv.\ in Appl.\ Math.\ \textbf{3}, no.\ 4, (\href{https://mathscinet.ams.org/mathscinet-getitem?mr=682631}{1982}), 463--490.
\bibitem{MP1990_SutoFH}
	A.~Süt\H{o}, \textit{Spectra of some almost periodic operators}, Number theory and physics (Les Houches, 1989), Springer Proc.\ Phys., \textbf{47}, (\href{https://mathscinet.ams.org/mathscinet-getitem?mr=1058459}{1990}), 162--169.
\bibitem{MP2000_TeschlReview}
	G.~Teschl, \textit{Jacobi operators and completely integrable nonlinear lattices}, Mathematical Surveys and Monographs \textbf{72}, (\href{https://mathscinet.ams.org/mathscinet-getitem?mr=1711536}{2000}), xvii+351.
\bibitem{TOP1970_Willard}
	S.~Willard, \textit{General topology}, Addison-Wesley Publishing Co., Reading, Mass.-London-Don Mills, Ont., (\href{https://mathscinet.ams.org/mathscinet-getitem?mr=264581}{1970}), xii+369 pp.
\bibitem{RA2014_Yeh}
	J.~Yeh, \textit{Real analysis:\ Theory of measure and integration}, Third edition, World Scientific Publishing Co.\ Pte.\ Ltd., Hackensack, NJ, (\href{https://mathscinet.ams.org/mathscinet-getitem?mr=3308472}{2014}), xxiv+815 pp.
\bibitem{MP2020_ZhaoAMO}
	X.~Zhao, \textit{Continuity of the spectrum of quasi-periodic Schrödinger operators with finitely differentiable potentials}, Ergodic Theory Dynam.\ Systems \textbf{40}, no.\ 2, (\href{https://mathscinet.ams.org/mathscinet-getitem?mr=4048305}{2020}), 564--576.
\end{thebibliography}
\end{document}